\documentclass[a4paper]{amsart}
\usepackage{graphicx}
\usepackage{epstopdf}

\usepackage{caption} 
\usepackage{subcaption} 
\usepackage{enumerate} 
\usepackage{fancyhdr} 

\usepackage{hyperref} 
\hypersetup{colorlinks, citecolor=black, urlcolor=blue,
         linkcolor=black, breaklinks=true, hypertexnames=true
         }

\usepackage{mathtools, commath}        

\newcommand{\NN}{{\mathbb{N}}}

\newcommand{\RR}{{\mathbb{R}}}

\newcommand{\EU}{{\mathbb{S}}}

\newcommand{\dpt}{\displaystyle}
\newcommand{\VC}{\mathcal{C}_0}

\newcommand{\p}{\mathfrak{p}}
\DeclareMathOperator{\Tr}{tr} 
\usepackage{xcolor} 

\newcommand{\citet}{\cite}
\newcommand{\citeauthor}{\cite}
\newcommand{\citeyear}{\cite}

\newcommand{\lessapprox}{\overset{<}{\approx}}
\newcommand{\gtrsim}{\overset{>}{\sim}}

\newtheorem{theorem}{Theorem}

\theoremstyle{definition}
\newtheorem{example}{Example}[section]

\newtheorem{definition}{Definition}[section]

\begin{document}

%

\markboth{B. F. F. Gon\c{c}alves, I. S. Labouriau and A. A. P. Rodrigues}{Bifurcations and canards in the \emph{FitzHugh-Nagumo} system: a {\isl tutorial}}

\title[Bifurcations and canards in the \emph{FitzHugh-Nagumo} system]{Bifurcations and canards in the {FitzHugh-Nagumo} system: \\  a tutorial in fast-slow dynamics}

\author[B.F.F. Gon\c{c}alves]{Bruno F. F. Gon\c{c}alves}

\address{B.F.F. Gon\c{c}alves --- Centro de Matem\'atica da Universidade do Porto\\ Rua do Campo Alegre 687, Porto 4169-007, Portugal}
\email{brunoffg9@gmail.com}

\author[I.S. Labouriau]{ Isabel S. Labouriau}

\address{I.S. Labouriau --- Centro de Matem\'atica da Universidade do Porto \\ Rua do Campo Alegre 687, Porto 4169-007, Portugal}
\email{islabour@fc.up.pt}

\author[A.A.P. Rodrigues]{Alexandre A. P. Rodrigues}
\address{A.A.P. Rodrigues --- Lisbon School of Economics \& Management\\Centro de Matem\'atica Aplicada \`a Previs\~{a}o e Decis\~{a}o Econ\'omica \\Rua do Quelhas 6,  Lisboa 1200-781, Portugal }
\email{
arodrigues@iseg.ulisboa.pt}

\begin{abstract}
    In this article, we study the FitzHugh-Nagumo $(1,1)$--fast-slow system where the vector fields associated to the fast/slow equations come from the reduction of the Hodgin-Huxley model for the nerve impulse. After deriving dynamical properties of the singular and regular cases, we perform a bifurcation analysis and we investigate how the  parameters (of the affine slow equation) impact the dynamics of the system. The study of codimension one bifurcations and the numerical locus of \emph{canards} concludes this case-study. All theoretical results are numerically illustrated.
\end{abstract}



\maketitle

\section{Introduction}\label{intro}
    Many animals have cells in their various physiological systems, such as the nervous, muscular, and cardiac systems, that are sensitive to certain electrical stimuli.
    These cells remain \emph{at rest} most of the time, react to electrical stimulation at a given moment, and then return to a resting state until they are stimulated again. 
    This oscillatory and excitability dynamics is present in neuronal activity and cardiac rhythm -- see for instance \citet{keener}.

  In  1952, Hodgkin \& Huxley  \citeyear{Hodgkin5}    
    developed a mathematical model that described the propagation of electrical signals along a squid's giant axon, known as the {\em Hodgkin-Huxley} (HH) model.
    The authors related the excitability of squid nerve cells and the resulting electrical impulse to the potential difference between the inside and outside of the squid's axon membrane  arising from the movement of sodium and potassium ions across the membrane.
    Using experimental data, they established the propagation of electrical signals along the squid's giant axon through a nonlinear system of four differential equations \cite{Hodgkin5}.

    In 1960--1961  Richard FitzHugh \cite{FitzHugh1,FitzHugh2}
     created a simplified version of the HH model, considering only two variables. 
     His goal was to simplify the HH model in order to make the dynamics of excitability more perceptible. The system developed by  by FitzHugh in 1961
      has the following form:
        \begin{equation}\label{Fitz}
            \begin{aligned}
                \frac{dx}{dt} &= a\left(x - \frac{x^3}{3} - y + I\right) \\
                \frac{dy}{dt} &= \frac{1}{a}\left(x - by - c\right),
            \end{aligned}
        \end{equation}
    %
    %
    \smallskip
    where $a\in \RR\backslash \{0\}, b,c\in \RR$ are constants and $I$ may be seen as the intensity of a stimulus applied to the axon membrane. Variables $x$ and $y$ of \eqref{Fitz} describe the fast evolution of the membrane voltage of a neuron as well as the activation of sodium channels, and the inactivation/activation of the sodium/potassium channels, respectively. They match with the pairs of variables $\left(v,m\right)$ and $\left(h,n\right)$ from the HH model, respectively -- cf. \cite{FitzHugh2,keener}. 
    
    During the sixties, Jin-Ichi Nagumo constructed an electrical circuit that recreated FitzHugh's system \cite{Nagumo}.
    Nowadays, due to the contributions of both researchers, this model is known as the FitzHugh-Nagumo (FH-N) model. This model has been extensively studied see, for instance, \citet{lacasa24} and references therein.
    
    By changing the time scale of system \eqref{Fitz}, considering $\tau=t/ a$ and $\varepsilon=1/a^2$, we obtain an equivalent version of the FH-N model as follows:
    \begin{equation}\label{epsFitz}
            \begin{aligned}
                \varepsilon\frac{dx}{d\tau} &= x - \frac{x^3}{3} - y + I \\
                \frac{dy}{d\tau} &= x - by - c.
            \end{aligned}
    \end{equation}
    For $0<\varepsilon\ll1$, or equivalently, $a\gg1$, system \eqref{epsFitz} is what we call a fast-slow system. 
    
    A \emph{fast-slow} system is a system of differential equations in which some variables have their derivatives with larger magnitude than others. This leads to a system with multiple time scales. The general approach to this type of systems starts by grouping the variables in two disjoint sets: fast variables and slow variables. This separation is introduced in system \eqref{epsFitz} by the parameter $\varepsilon$.

    In Doss-Bachelet {\em et al.} \citet{artBachelet}    
    the authors notice that the FH-N model presents two distinct dynamics near a Hopf bifurcation: an attracting focus where all trajectories converge to it if the parameter is less than the bifurcation value; and an attracting relaxation-oscillation periodic solution (limit cycle) if the parameter is greater than the bifurcation point. The work of \citeauthor{artBachelet} proceeds with the analysis of the dynamics and applies it to generate bursting oscillations from two coupled FH-N systems, thus leaving the FH-N model little explored.

    \subsection{Literature}
        In this section, we review the FH-N equation, as well as its derivations, as one paradigmatic example of mathematical modelling.
        Although it has been motivated by neuroscience, the FH-N model exists in many variations of the original system and has been studied in several fields due to its simplicity and rich dynamics. 
        Motivated by its broad applicability and generic dynamics, the FH-N system has triggered a huge body of research works characterising their bifurcations and the associated dynamical regimes. 
        
        There is an abundance of references in the literature. In what follows, we choose to mention only a few for clarity and the choice is  based on our personal preferences. The reader interested in further detail and/or more examples can use the references within those we mention specially those in the review \citet{lacasa24}.
        
        At the first glance, the simplifications introduced in the passage from the sophisticated HH model to the FH-N model might have a price. Phillipson and Schuster \cite{phillipson} explore numeric and analytically the four-dimensional Hodgkin–Huxley model and demonstrate that, under precise conditions, the FH-N equations can provide quantitative predictions in close agreement with the classic HH equations. 
        
        The bifurcations of FH-N system have been extensively analysed as a fast-slow system
in \cite{chow, canards_numerics, MMO_Krupa,  FastCoupling_Kristian, kuehn,  MMO_Desroches, FEN, KS01a} uncovering complex patterns of behaviours that emerge from the simple and/or coupled FH-N model. They reveal how changes in the parameters can lead to different dynamical regimes, including equilibria, limit cycles, \emph{canards} and chaos. Stability techniques allow the identification of critical points where small perturbations may lead to qualitatively different behaviours, shedding light on the robustness and sensitivity of excitable systems.  Systems exhibiting excitability often display a threshold behaviour, where a certain level of input of perturbation is required to trigger the response -- cf. \cite{artBachelet}.
        
        We would like to refer the paper \cite{model2003global} for the study, not using the fast-slow approach,  of \eqref{Fitz}, where global bifurcations are presented in the context of the global bifurcation diagram. Codimension one and two bifurcations of Hopf, homoclinic, saddle-node and Bogdanov-Takens are obtained.   
        
        The work \cite{abdelouahab2019complex} investigates the complex phenomena of \emph{canard} explosion with mixed-mode oscillations, which can be observed in a fractional-order FH-N model.
        The study has highlighted the appearance of patterns of solutions with an increasing number of small-amplitude oscillations in each of such patterns, as one parameter of the fractional-order FH-N model is varied.
        
        The paper \cite{llibre2015periodic} extends the FH-N model by including a time-varying threshold between electrical ``silence'' and electrical firing, as well as a periodic forcing for the intensity of the stimulus. Authors provide sufficient conditions for the existence of periodic solutions in such differential systems. This may be related with the work \cite{wang2007dynamical}, where the author explores the dynamical behaviour of non-autonomous, almost-periodic discrete FH-N system defined on infinite lattices. The non-autonomous infinite-dimensional system has a uniform attractor which attracts all solutions uniformly with respect to the translations of external terms. Authors establish the upper semi-continuity of uniform attractors when the infinite-dimensional system is approached by a family of finite-dimensional systems. 
        
        In \cite{dong2015identification} authors propose a new method for the identification of the FH-N model dynamics via deterministic learning (with diffusion). The FH-N model is then described by partial differential equations which are transformed into a set of ordinary differential equations whose dynamics has been studied. Authors found spiral waves and an accurate identification of dynamics underlying the recurrent trajectory corresponding to any spatial point is achieved. 
        
        Spiking and bursting  patterns observed in nerve membranes seem to be important when we investigate information representation model in the brain. Many topologically different bursting responses are observed in the mathematical models and their related bifurcation mechanisms have been clarified in \cite{tsuji2004design}. In this article, the authors propose a design method to generate bursting responses in FitzHugh–Nagumo model with a simple periodic external force based on bifurcation analysis. Some effective parameter perturbations for the amplitude of the external input are given from the two-parameter bifurcation diagram.

    \subsection{Novelty and structure of the article}
        Our goal with the present article is to exhibit the various dynamics in the multiple time scale model and motivate their existence in the light of methods from geometric singular perturbation theory and bifurcation theory.  
        This article studies the dynamics of a fast-slow system derived from the FH-N model according to the geometric singular perturbation theory and the bifurcation theory methods. We provide an analytic proof that the fast-slow FH-N system presents relaxation oscillation dynamics as well as periodic solutions induced by Hopf bifurcation -- we fully analyse the type of bifurcations and we check rigorously all conditions for the emergence of a homoclinic orbit and \emph{canards} connecting these two distinct dynamics. We illustrate each result with numerical computations.

        This article is organised as follows. In Section \ref{sec:sing} we introduce general concepts of fast-slow systems and geometric singular perturbation theory and apply them to the analysis of the singular case of the FH-N system. A relaxation-oscillation periodic solution is studied with an estimate of its period in Section \ref{sec:per}. In Section \ref{sec:reg} we describe Fenichel's theorem and link singular to regular dynamics and the way the latter converges to the former as $\varepsilon\to0$ applying it to FH-N. In Section \ref{sec:bif} we perform a bifurcation analysis of FH-N and investigate the influence of parameters in the qualitatively different dynamics of the model. The study of singularities and the emergence of \emph{canard} trajectories in Section \ref{sec:canards} concludes this tutorial. Section \ref{discussion} discusses the results presented and states the natural continuation of the present work.
        
         We have endeavoured to make a self contained exposition bringing together all topics related to the proofs. We have drawn illustrative figures to make the paper easily readable. All figures in this article were created through numerical simulations conducted in \emph{Matlab}, using integration functions such as \emph{ode15s} or \emph{ode23s}, and, in more sensitive examples, the \emph{MatCont} toolbox developed by Willy Govaerts, Yuri A. Kuznetsov and Hil G.E. Meijer \cite{matcont}.

\section{Multiple Time Scale Theory and the FitzHugh-Nagumo System}\label{ch1}
    In this section, general concepts of two time scale systems are introduced and illustrated with the FH-N system. Further results on \emph{fast-slow} systems may be found in \citet{kuehn}. We start by introducing some necessary definitions.
    \begin{definition}\label{def:2.1}
         For $m, n\in \NN$, a $(m,n)$--fast-slow system is a system of $m + n$ ordinary differential equations of the form:
         \begin{equation}\label{deflento}
             \begin{aligned}
                 &\varepsilon\,\dfrac{\mathrm{d}x}{\mathrm{d}\tau}=\varepsilon\Dot{x}=f(x,y,\varepsilon),\\
                 &\phantom{\varepsilon\,}\dfrac{\mathrm{d}y}{\mathrm{d}\tau}=\phantom{\varepsilon}\Dot{y}=g(x,y,\varepsilon),
             \end{aligned}
         \end{equation}
    \end{definition}
    \noindent where $f:\RR^m\times\RR^n\times\RR\to\RR^m$, $g:\RR^m\times\RR^n\times\RR\to\RR^n$ are $C^2$ and $0<\varepsilon\ll1.$
    We denote the \textbf{fast} variable(s) by $x \in \RR^m$ and $y \in \RR^n$ is called the \textbf{slow} variable(s). 
    Considering $t=\tau/\varepsilon$, system \eqref{deflento} is equivalent to:
        \begin{equation}\label{defrapido}
            \begin{aligned}
                &\dfrac{\mathrm{d}x}{\mathrm{d}t}=x'=\phantom{\varepsilon\,} f(x,y,\varepsilon),\\
                &\dfrac{\mathrm{d}y}{\mathrm{d}t}=y'=\varepsilon\, g(x,y,\varepsilon).
            \end{aligned}
        \end{equation}

    The \textbf{fast time} scale is represented by $t$ while $\tau$ refers to the \textbf{slow time} scale. In order to better understand the concepts and the terminology associated with these systems, the following definitions and results are illustrated by the FH-N (1,1)--fast-slow system.
    \begin{example}\label{ex:fhn}
        As described in \citet{artBachelet}, consider the FH-N (1,1)--fast-slow system:
        \begin{equation}
            \begin{aligned}
            \label{fhn}
                &\text{Slow time scale } \left(\tau\right) \qquad\qquad && \text{Fast time scale } {\textstyle\left(t=\tau/\varepsilon\right)} \\
                &\varepsilon\,\dfrac{\mathrm{d}x}{\mathrm{d}\tau}=\varepsilon\Dot{x}=f(x,y,\varepsilon) \qquad\qquad &&\dfrac{\mathrm{d}x}{\mathrm{d}t}=x'=\phantom{\varepsilon}\,f(x,y,\varepsilon)\\
                &\phantom{\varepsilon\,}\dfrac{\mathrm{d}y}{\mathrm{d}\tau}=\phantom{\varepsilon}\Dot{y}=g(x,y,\varepsilon) \qquad\qquad &&\dfrac{\mathrm{d}y}{\mathrm{d}t}=y'=\varepsilon\,g(x,y,\varepsilon),
            \end{aligned}
        \end{equation}
        \noindent where $f,g:\; \RR^3\longrightarrow\RR$ are given by 
        \begin{equation*}
            f(x,y,\varepsilon)=-y+4x-x^3, \quad g(x,y,\varepsilon)=x-by-c    
        \end{equation*}
        with parameters $b,c \in{\RR}$ and $0<\varepsilon\ll1$. $\diamond$ 
    \end{example}

    Both cubic and affine functions of Example \ref{ex:fhn}, $f$ and $g$, come naturally from the reduction of the HH model \cite{keener}. Note that in FH-N the functions $f$ and $g$ do not depend on $\varepsilon$. System \eqref{fhn} will be referred to as the FH-N system in two variables.
\section{The Singular Case }\label{sec:sing}

    When analysing a fast-slow system, it is useful to consider the singular case $\varepsilon=0$. Techniques such as nullclines analysis provide a geometric understanding of the model's behaviour.
    \begin{definition}
        For $\varepsilon=0$, system \eqref{deflento} in the slow time scale becomes a system of differential-algebraic equations:
        \begin{equation}\label{deffaselenta}
            \begin{aligned}
                &0=f(x,y,\varepsilon),\\
                &\Dot{y}=g(x,y,\varepsilon).
            \end{aligned}
        \end{equation}
    \end{definition}
    System \eqref{deffaselenta} is called \textbf{slow subsystem}, its equations are referred to as the \textbf{reduced equations} and its flow the \textbf{slow/reduced flow}.
    \begin{definition}
        For $\varepsilon=0$, system \eqref{defrapido} in fast time is called \textbf{fast subsystem}:
        \begin{equation}\label{deffaserapida}
            \begin{aligned}
                 &x'=f(x,y,\varepsilon),\\
                 &y'=0.
             \end{aligned}
        \end{equation}
    \end{definition}
    \noindent The set of equations in system \eqref{deffaserapida} is referred to as the \textbf{layer equations} and its flow as the \textbf{fast flow}.

\begin{figure}
        \centering
        \begin{subfigure}{.5\linewidth}
            \centering
            \includegraphics[width=\linewidth]{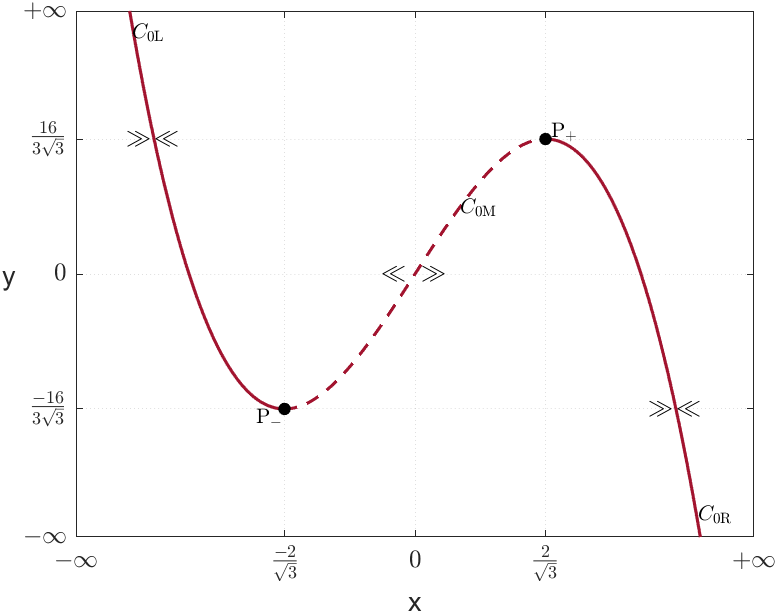}
            \caption{ \small }
            \label{fig:C0_2d}
        \end{subfigure}%
        \begin{subfigure}{.5\textwidth}
            \centering
            \includegraphics[width=\linewidth]{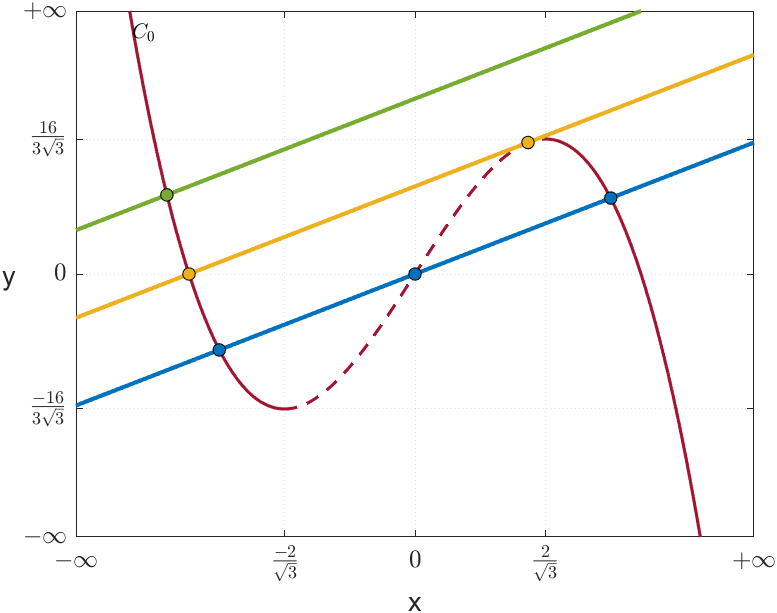}
            \caption{ \small }
            \label{fig:C0_2d_Eqs}
        \end{subfigure}
        \caption{ \small  (a) - Critical manifold of the system \eqref{fhn}: Attracting (solid red) and repelling (dashed red) regions; Fold points (black dots). (b) - Geometrical argument for the existence of at least one and at most three equilibrium points in the system \eqref{fhn}.}
        \label{fig:C0_2d_e_Eqs}
    \end{figure}

    Systems \eqref{deffaselenta} and \eqref{deffaserapida} allow the understanding of the dynamics of \eqref{deflento} and \eqref{defrapido} for small $\varepsilon$, $0<\varepsilon\ll1$, since the flow may be seen as either the fast or the slow flow, depending on the region in the phase space. Near the set defined by $f\left(x,y,0\right)=0$, one expects the flow to be \emph{$C^2$--close to} the slow flow as we will see later in Section \ref{sec:reg}.
    %
    %
    \begin{definition}
        For $m, n\in \NN$, given a $(m,n)$--fast-slow system as in \eqref{def:2.1}, the \textbf{critical manifold} is the set
        \begin{equation}{\label{c0}}
            \VC = \left\{(x,y)\in \RR^m\times\RR^n :f(x,y,0)=0\right\}.
        \end{equation}
    \end{definition}
    \bigbreak
    \begin{example}
        The critical manifold for the FH-N system \eqref{fhn} is $$\VC = \left\{(x,y)\in \RR^2 :f(x,y,0)=0\right\} = \left\{(x,y)\in \RR^2 :y=4x-x^3\right\},$$
        a cubic shaped curve, parametrised by
        \begin{equation}\label{eq:C0_param}
            x\mapsto (x,4x-x^3)=(x,y),
        \end{equation}
        as depicted in Figure (\ref{fig:C0_2d}). The map $\varphi(x)=4x-x^3$ is odd $ \left(\forall x\in \RR,-\varphi\left(x\right)=\varphi\left(-x\right)\right)$.
        $\diamond$
    \end{example}

    All points in $\VC$ are equilibria of the fast subsystem. Through \emph{Hartman-Grobman}'s Theorem \cite{hartman}, we may define the (normal) hyperbolicity and stability of the points in $\VC$ with respect to the fast variables, in the following way:
   
    \begin{definition}
        For $m, n\in \NN$, given a $(m,n)$--fast-slow system with critical manifold $\VC$:
        \begin{enumerate}
            \item The subset $S\subset\VC$ is said to be \textbf{normally hyperbolic} if for all $\p\in{S}$, we have 
            \begin{equation*}
                \left(\mathsf{D}_xf\right)\left(\p,0\right) \text{ has no eigenvalues with zero real part,}
            \end{equation*}
            where $\left(\mathsf{D}_xf\right)$ denotes the Jacobian matrix of $f$ with respect to $x$.

            \medskip
        
            \item A normally hyperbolic set $S$ is said to be \textbf{attracting}/\textbf{repelling} if, for all $\p\in{S}$, all eigenvalues of $\left(\mathsf{D}_xf\right)\left(\p,0\right)$ have negative/positive real part, respectively. It is of \textbf{saddle type} if $\left(\mathsf{D}_xf\right)\left(\p,0\right)$ has eigenvalues with   both positive and negative real part.
        \end{enumerate}
    \end{definition}
    
    %
    %
    When considering a (1,1)--fast-slow system, since there is only one fast variable,  we have $\mathsf{D}_xf = \dfrac{\partial f}{\partial x}$. For convenience, we use the notation $\varphi_\sigma$ when referring to the partial derivative of a real map $\varphi$ in order to its variable $\sigma$.
    \begin{definition}\label{def:dobra}
        Given a (1,1)--fast-slow system with critical manifold $\VC$, the point $\p\in\VC$ is said to be a \textbf{fold point} if:
        \begin{equation*}
            f_x\left(\p,0\right)=0\,,\quad f_{xx}\left(\p,0\right)\neq0 \quad\text{ and }\quad f_y\left(\p,0\right)\neq0.
        \end{equation*}
        If $g\left(\p,0\right)\neq0$, the fold point is called \textbf{regular}.
    \end{definition}
    \begin{example}\label{casos}
        With respect to \eqref{fhn}, we have $f\left(x,y,0\right)=-y+4x-x^3$ and, for $\p\in{\VC}$, one has: $$f_x\left(\p,0\right)=4-3x^2\neq0\qquad\text{ if }\,\, x\neq\pm\dfrac{2}{\sqrt{3}}.$$
        Therefore, the equilibria $P_+ =\left( 2/\sqrt{3},\,16/3\sqrt{3}\right)$ and the set $P_- =\left( -2/\sqrt{3},\,-16/3\sqrt{3}\right)$ are fold points and $S_0=\VC\setminus\left\{\left(x,\varphi\left(x\right)\right):\,x=\pm 2/\sqrt{3}\right\}$ is normally hyperbolic.
        We may split $S_0$ into three distinct open subsets:
        \begin{equation*}
            \begin{aligned}
                \mathcal{C}_{0\text{L}}&= \VC \cap {\textstyle \left\{\left(x,y\right)\in\RR^2: x<\dfrac{-2}{\sqrt{3}}  \right\}}, \\
                \mathcal{C}_{0\text{M}}&= \VC \cap {\textstyle \left\{\left(x,y\right)\in\RR^2: \dfrac{-2}{\sqrt{3}}<x<\dfrac{2}{\sqrt{3}}  \right\}}\quad \text{and} \\
                \mathcal{C}_{0\text{R}}&= \VC \cap {\textstyle \left\{\left(x,y\right)\in\RR^2: x>\dfrac{2}{\sqrt{3}}  \right\}}. \\
            \end{aligned}
        \end{equation*}

        Therefore, we get 
        \begin{equation}
           \label{def: S0}
    S_0=\mathcal{C}_{0\text{L}}\cup\mathcal{C}_{0\text{M}}\cup\mathcal{C}_{0\text{R}},    
        \end{equation}
        where $\mathcal{C}_{0\text{L}}$ and $\mathcal{C}_{0\text{R}}$ are attracting subsets and $\mathcal{C}_{0\text{M}}$ is repelling. This means that in the fast flow, trajectories of \eqref{fhn} move towards either $\mathcal{C}_{0\text{L}}$ or $\mathcal{C}_{0\text{R}}$.
        $\diamond$ 
    \end{example}
        %

    Since the critical manifold $\VC$ is parametrised by the fast variable $x$ (cf. Eq.\eqref{eq:C0_param}), it is possible to explicitly define the slow flow through the fast variable $x$. This allows the stability analysis of the fast and slow flows expressed in terms of the fast variable. \\  

    In system \eqref{fhn}, by differentiating $f(x,y,0)=-y+4x-x^3=0$ with respect to $\tau$ we have: $$\dfrac{\partial f}{\partial x} \Dot{x} + \dfrac{\partial f}{\partial y} \Dot{y} = (4-3x^2)\Dot{x} - \Dot{y}=0.$$
     
    Since $\Dot{y}=x-by-c$, then $$\Dot{x}=\dfrac{x-by-c}{(4-3x^2)}.$$ 
        
    Taking into account that $y=4x-x^3$, we may conclude that 
    \begin{equation}\label{expVC}
        \Dot{x} = \dfrac{bx^3 + (1-4b)x + c}{4-3x^2} \eqqcolon \psi\left(x,b,c\right).
    \end{equation}
        
    As expected, this explicit form of the slow subsystem is not defined at the fold points $x=\pm 2/\sqrt{3}$.
        
    If the fold points are regular (cf. Def. \ref{def:dobra}), then all equilibria $\left(x^*,y^*\right)\in\VC$ of the system \eqref{fhn} belong necessarily to $S_0$. Their stability in the slow direction may be determined in the following way: 
    \begin{equation*}
        \dfrac{\partial \psi}{\partial x}= \psi_x\left(x,b,c\right)=\dfrac{1-b\left(4-3x^2\right)+6x\,\psi\left(x,b,c\right)}{\left(4-3x^2\right)}.
    \end{equation*}

    Since $\psi\left(x^*,b,c\right)=0$, then
        
    \begin{equation*}
        \psi_x \left(x^*,b,c\right)=\dfrac{1-b\left(4-3{x^*}^2\right)}{4-3{x^*}^2}\,.
    \end{equation*}
        
    For $b>0$, if $|x^*|>2/\sqrt{3}$,  the equilibrium attracts the slow flow. Otherwise, if $|x^*|<2/\sqrt{3}$ the equilibrium may either attract or repel the slow flow depending on the value of $1-b\left(4-3{x^*}^2\right)$ being  negative or positive, respectively. Thus, the equilibrium point $\left(x^*,y^*\right)$ of \eqref{fhn} is stable if $\left(x^*,y^*\right)\in \mathcal{C}_{0\text{L}}\cup\mathcal{C}_{0\text{R}}$ and unstable if $\left(x^*,y^*\right)\in\mathcal{C}_{0\text{M}}$.
        
    For $b<0$, the stability in the slow direction is the opposite of the previous one. Hence, the equilibrium point $\left(x^*,y^*\right)$ is unstable if $\left(x^*,y^*\right)\in\mathcal{C}_{0\text{M}}$ and  either stable or unstable of saddle type if $\left(x^*,y^*\right)\in \mathcal{C}_{0\text{L}}\cup\mathcal{C}_{0\text{R}}$.
       
    In Example \ref{casos}, we have disregarded the case where a fold point is also an equilibrium of the system, \emph{i.e.}, $x^*=\pm 2/\sqrt{3}$. This particular situation requires special attention and will be analysed later in Section \ref{sec:canards}. For the remainder of this section, we   describe the dynamics in the cases where all equilibria  {lie in the set $S_0$ defined in \eqref{def: S0}.}
    %
\subsection{Dynamics}
    Depending on the values of $b$ and $c$, there   exists at least one and at most three equilibria, as result of the intersection of the cubic curve, $f\left(x,y,0\right)=0$, with the line, $\dot{y}=0$ (cf. Figure \ref{fig:C0_2d_Eqs}).
    Moreover, when the system has three equilibria only the following scenarios may occur: either all three equilibria belong to $\mathcal{C}_{0\text{M}}$ or else they are each one in a distinct region, $\mathcal{C}_{0\text{L}}$, $\mathcal{C}_{0\text{M}}$ and $\mathcal{C}_{0\text{R}}$. \\
     
\begin{figure}
        \centering
        \begin{subfigure}{0.3\linewidth}
            \includegraphics[width=\linewidth]{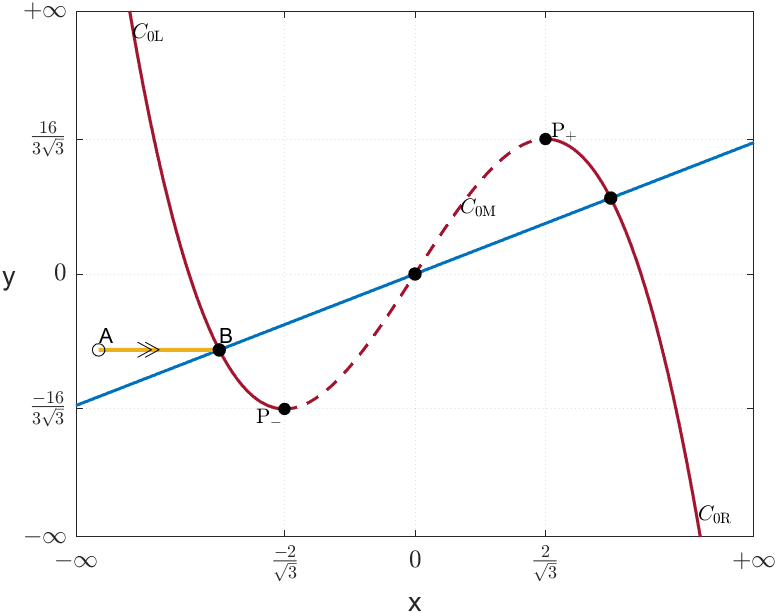}
            \caption{ \small  Case (i).}
            \label{fig:caso1}
        \end{subfigure}
        \begin{subfigure}{0.3\linewidth}
            \includegraphics[width=\linewidth]{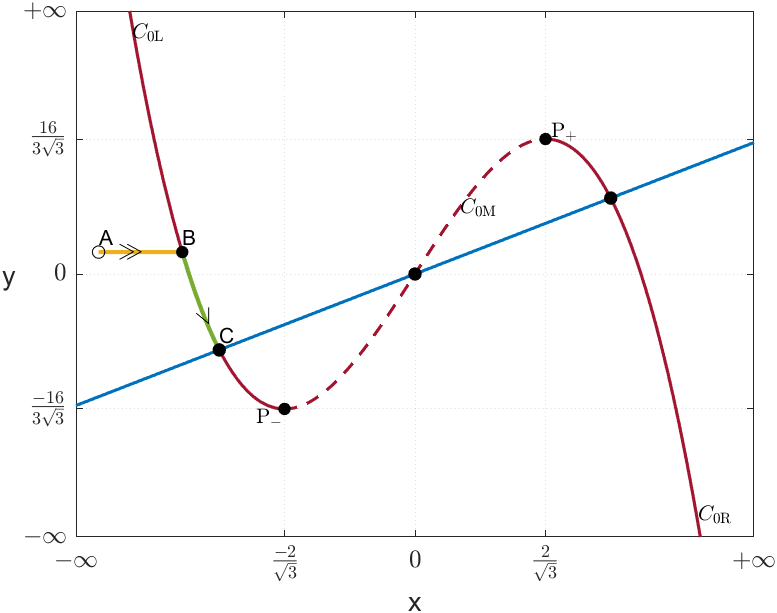}
            \caption{ \small  Case (ii).}
            \label{fig:caso2}
        \end{subfigure}
        \begin{subfigure}{0.3\linewidth}
            \includegraphics[width=\linewidth]{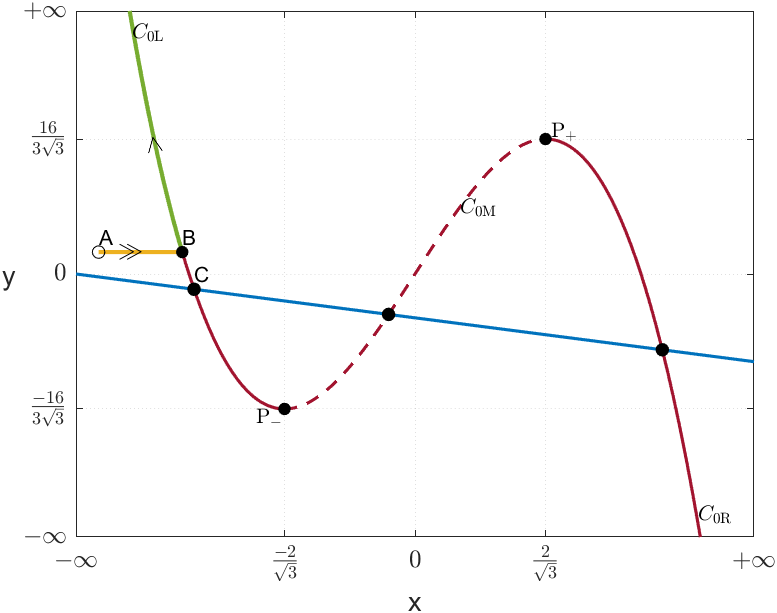}
    		\caption{ \small  Case (iii) - $y_a>y_c$.}
    		\label{fig:caso3a}
        \end{subfigure}
	\vfill
        \begin{subfigure}{0.3\linewidth}
            \includegraphics[width=\linewidth]{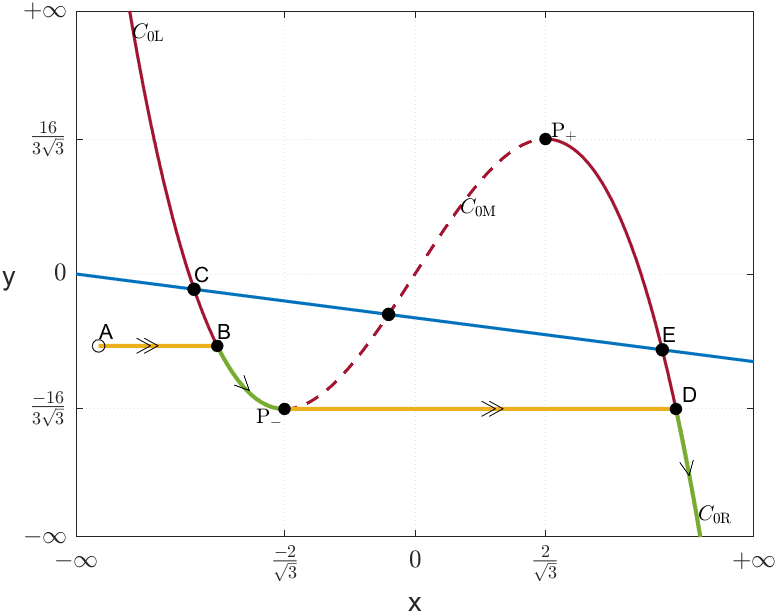}
            \caption{ \small  Case (iii) - $y_a<y_c$.}
            \label{fig:caso3b}
        \end{subfigure}
        \begin{subfigure}{0.3\linewidth}
            \includegraphics[width=\linewidth]{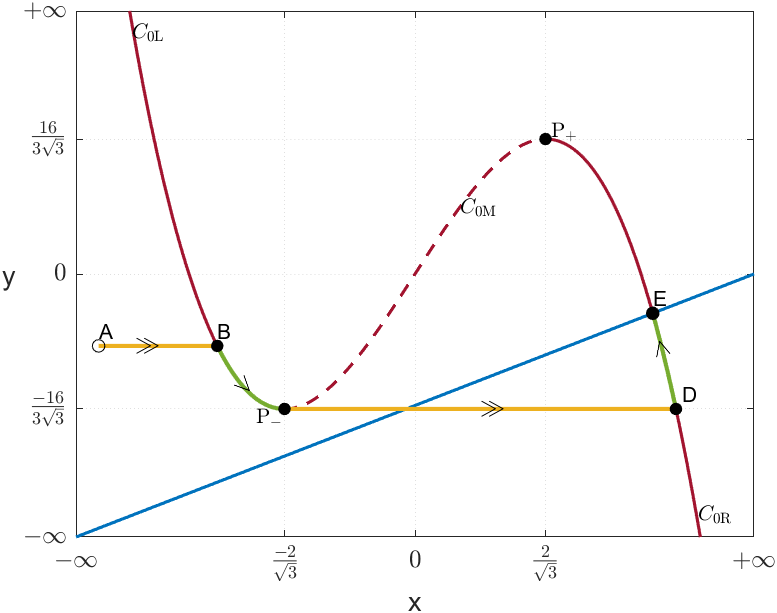}
            \caption{ \small  Case (iv) - $E$ stable.}
            \label{fig:caso4es}
        \end{subfigure}
        \begin{subfigure}{0.3\linewidth}
            \includegraphics[width=\linewidth]{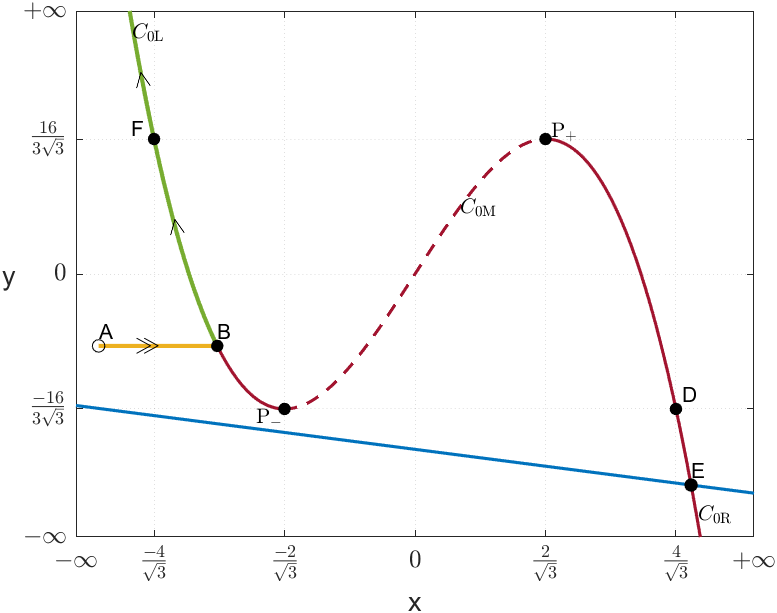}
            \caption{ \small  Case (iv) - $E$ unstable.}
            \label{fig:caso4ins}
        \end{subfigure}
        \caption{ \small  Examples of Cases (i)--(iv). Critical manifold in red, attracting regions in solid line and repelling region in dashed line. Line $\dot{y}=0$ in blue. Fast trajectories in yellow and slow trajectories in green. Points of interest/equilibrium/fold in black. Initial point, $A$, outlined in black.}
        \label{fig:casos_1-4}
    \end{figure}

    Fix $A= \left(x_a,y_a\right)\in\RR^2\setminus\VC$.   Consider the half trajectory $\gamma_A$  for $t \geq 0$ that starts at $A$. Since $y'=0$, the fast flow runs \emph{horizontally} towards $\mathcal{C}_{0\text{L}}$ or $\mathcal{C}_{0\text{R}}$.  Without loss of generality, suppose that $\gamma_A$ moves towards $\mathcal{C}_{0\text{L}}$ and let $B=\left(x_b,y_a\right)\in\mathcal{C}_{0\text{L}}$ be the intersection point of $\gamma_A$ with $\VC$. Here, five different scenarios may occur (cf. Figure \ref{fig:casos_1-4}):
    \begin{enumerate}[\hspace{0.5cm} (i)]
        \item The point $B$ is an equilibrium;
        \item The point $B$ is not an equilibrium, but there exists a stable equilibrium $C=\left(x_c,y_c\right)\in\mathcal{C}_{0\text{L}}\setminus {B}$;
        \item The point $B$ is not an equilibrium, but there exists an unstable equilibrium $C=\left(x_c,y_c\right)\in\mathcal{C}_{0\text{L}}\setminus {B}$;
        \item There is no  equilibrium in $\mathcal{C}_{0\text{L}}$, but it exists in $\mathcal{C}_{0\text{R}}$;
        \item All equilibria lie in $\mathcal{C}_{0\text{M}}$.
    \end{enumerate}
   
    Scenario (i) is trivial. If $B=\left(x_b,y_a\right)$ is an equilibrium of \eqref{fhn}, then by definition $\left(\Dot{x},\Dot{y}\right)=0$, and therefore, the trajectory initiated at $A$  accumulates at $B$ (cf. Figure \ref{fig:caso1}).\\
    
    Scenario (ii) is not much different. Upon reaching the point $B$, the trajectory $\gamma_A$ moves through the slow flow of the system defined by $\Dot{y}=x-by-c$. Since the equilibrium point $C$ is attracting, and since both points $B$ and $C$ are in $\mathcal{C}_{0\text{L}}$, the trajectory $\gamma_A$ moves from $B$ to $C$ along the branch $\mathcal{C}_{0\text{L}}$ of the critical manifold, where it accumulates  (cf. Figure \ref{fig:caso2}).\\
  
    Scenario (iii):  First note that $\dfrac{\partial g}{\partial y}\left(x,y,0\right)=-b>0$ is constant for all $(x,y)\in \VC$. Since $\mathcal{C}_{0\text{L}}$ and $\mathcal{C}_{0\text{R}}$ attract the slow flow we have that any equilibrium point in either $\mathcal{C}_{0\text{L}}$ or $\mathcal{C}_{0\text{R}}$ is a saddle. In particular, the equilibrium point $C$ is a saddle.
   
    Since $\mathcal{C}_{0\text{L}}$ attracts the fast flow, we know that $\gamma_A$ follows $\mathcal{C}_{0\text{L}}$ away from $C$. Thus, if $y_a>y_c$, then $\gamma_A$ remains in $\mathcal{C}_{0\text{L}}$ and the second coordinate of $\gamma_A\left(\tau\right)$ tends to $+\infty$ as in Figure \ref{fig:caso3a}.  
    
    However, if $y_a<y_c$, the trajectory follows the slow flow to the fold point $P_-=\left(-2/\sqrt{3},\,-16/3\sqrt{3}\right)$. Here, since the slow flow is not defined at the fold points (cf. \eqref{expVC}), the trajectory $\gamma_A$ moves through the fast flow, \emph{horizontally}, to the point $D=\left(4/\sqrt{3},\,-16/3\sqrt{3}\right)\in\mathcal{C}_{0\text{R}}$ as in Figure \ref{fig:caso3b}. In this case there are two possibilities.
     
    The first possibility is that there is an equilibrium point $E=\left(x_e,y_e\right)\in\mathcal{C}_{0\text{R}}$ and that both $y_c$ and $y_e$ lie outside the interval $\left(-16/3\sqrt{3},\,\,16/3\sqrt{3}\right)$ delimited by the second coordinates of the two fold points. Then $\gamma_A$ follows the slow flow up to the other fold point $P_+=\left(2/\sqrt{3},\,\, 16/3\sqrt{3}\right)$, where it jumps back to $\mathcal{C}_{0\text{L}}$ and continues cycling around between the two branches. The trajectory $\gamma_A$ presents dynamics similar to case (v).
     
    In any other case, i.e., if either there is no equilibrium in $\mathcal{C}_{0\text{R}}$ or if one of the equilibria has its second coordinate in the interval $\left(-16/3\sqrt{3},\,\,16/3\sqrt{3}\right)$, then the second coordinate of $\gamma_A\left(\tau\right)$ tends to $\pm\infty$ as $\tau$ tends to $+\infty$. This is because, if $E=\left(x_e,y_e\right)\in\mathcal{C}_{0\text{R}}$ with $y_e<-16/3\sqrt{3}$ then the trajectory follows $\mathcal{C}_{0\text{R}}$ up to the fold point $P_+$, jumps to $\mathcal{C}_{0\text{L}}$ and follows it up, away from $C$ and with its second coordinate going to $+\infty$. If either there is no equilibrium in $\mathcal{C}_{0\text{R}}$ or if $E=\left(x_e,y_e\right)\in\mathcal{C}_{0\text{R}}$ with $y_e>-16/3\sqrt{3}$, then $\gamma_A$ follows $\mathcal{C}_{0\text{R}}$ down with its second coordinate going to $-\infty$ as in Figure \ref{fig:caso3b}.\\
     
    Scenario (iv): let $E=\left(x_e,y_e\right)$ be an equilibrium point in $\mathcal{C}_{0\text{R}}$. Geometrically, the point $E$ is the intersection of the cubic curve $x'=0$ with the line $\Dot{y}=0$. Since there are no more equilibria in $\mathcal{C}_{0\text{L}}$, we may observe that, if $E$ is unstable, then $y_a>y_e$.
    The proof is simple: if $E$ is unstable, then $\psi_x<0$  where $\psi$ was defined in \eqref{expVC}. Since $f_x\left(E,0\right)>0$, it follows that $b<0$ and therefore the line has a negative slope. Now, assume by contradiction that $y_b<y_e$. The line $\Dot{y}=0$ necessarily intersects the cubic surface in $\mathcal{C}_{0\text{L}}$ and hence, there exists an equilibrium in $\mathcal{C}_{0\text{L}}$ giving our contradiction.
    Thus, at $B$, we have $\Dot{y}>0$ and therefore, the trajectory $\gamma_A$ remains in $\mathcal{C}_{0\text{L}}$, with the first coordinate of $\gamma_A\left(\tau\right)$ going to $-\infty$ as $\tau\rightarrow\infty$ while the second coordinate of $\gamma_A\left(\tau\right)$ goes to $\infty$ { as in Figure \ref{fig:caso4ins}}. 
    If the point $E$ is stable, $\gamma_A$ travels along $\mathcal{C}_{0\text{L}}$ from $B$ to $P_-$, jumps in fast time to point $D$ in the branch $\mathcal{C}_{0\text{R}}$, ending at $E$ {(cf. Figure \ref{fig:caso4es})}. \\

\begin{figure}
        \centering
        \begin{subfigure}{.5\linewidth}
            \centering
            \includegraphics[width=\linewidth]{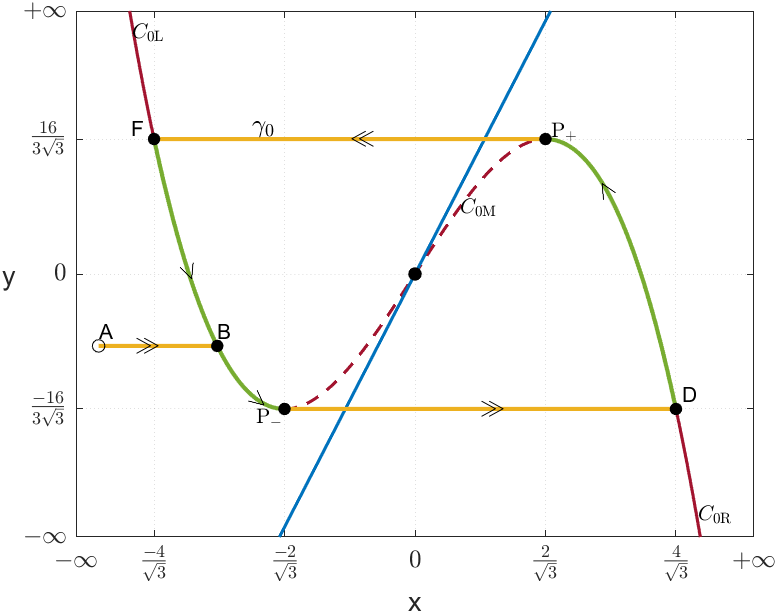}
            \caption{ \small  Case (v) - limit cycle.}
            \label{fig:caso5}
        \end{subfigure}%
        \begin{subfigure}{.5\textwidth}
            \centering
            \includegraphics[width=\linewidth]{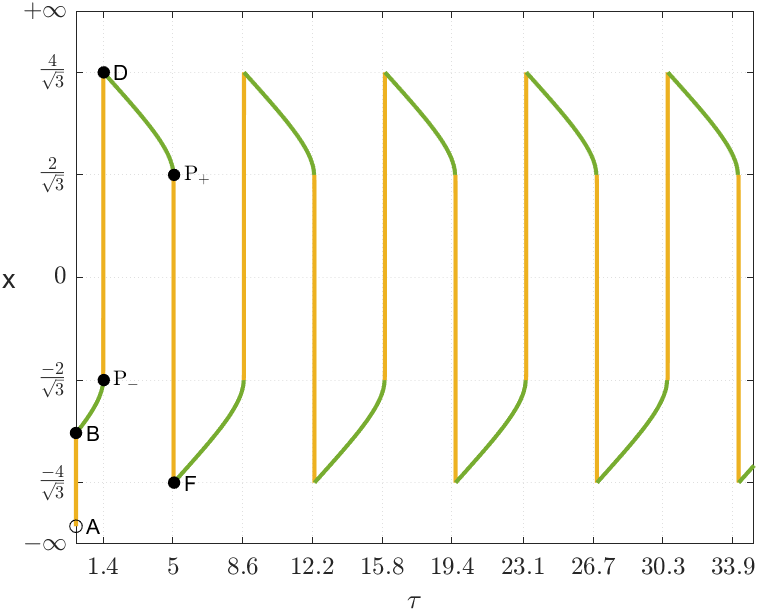}
            \caption{ \small  $x$ as function of $\tau$.}
            \label{fig:caso5_xtau}
        \end{subfigure}
        \caption{ \small  Example of case (v) with parameters $b=0.2$ and $c=0$. Periodic solution $\gamma_0$ composed of two slow trajectories (green) and two fast trajectories (yellow) with a period $\mathcal{T}_{\gamma_0}\approx 3.61$ time units ($\tau$).}
        \label{fig:caso5extau}
    \end{figure}

    Scenario (v): by using the symmetry of $f$ ($f$ is odd), one may check easily that the parameter $b$ is positive and that there are either three equilibria or only one in $\mathcal{C}_{0\text{M}}$. 
    In both cases, since the points belong to the unstable region of $\VC$, the dynamics of $\gamma_A$ does depend on the number of equilibria in the system. Therefore, we can observe that any point in the region $\mathcal{C}_{0\text{L}}$ satisfies $\Dot{y}<0$, so the trajectory $\gamma_A$ moves along $\mathcal{C}_{0\text{L}}$ from point $B$ to the fold point $P_-$, and then \emph{jumps}, in fast time, to point $D$.
    Similarly, $\gamma_A$ continues to $P_+$ and then \emph{jumps} to point $F$. Thus,  $\gamma_A$ goes through points $F,P_-,D$ and $P_+$ periodically and, therefore, $\gamma_A$ is a periodic solution.
    The trajectory neither diverges nor converges to an equilibrium point (cf. Figure \ref{fig:caso5}). 

    Any trajectory starting at any point in $\RR^2\setminus\VC$ moves in fast time to a point in one of the stable regions of the critical manifold $\VC$.
    In the previous scenarios, we have assumed that this point belongs to the branch $\mathcal{C}_{0\text{L}}$, but the dynamics for the case where the point is located on the branch $\mathcal{C}_{0\text{R}}$ is analogous to the former.
    %
\subsection{Period of the limit cycle}\label{sec:per}
    In Scenario (iii) when the equilibria of $\mathcal{C}_{0\text{L}}\cup\mathcal{C}_{0\text{R}}$ satisfy $|y|>\frac{16}{3\sqrt{3}}$ and in Scenario (v), we have observed the emergence of a periodic solution alternating between the slow and the fast flows.
    Now, in order to estimate its period, we are going to determine the points responsible for linking slow and fast flows.\\
      
    The transition from the slow flow to the fast one occurs at the fold points $ P_\pm=\pm\textstyle\left(2/\sqrt{3},\,16/3\sqrt{3}\right)$.  Similarly, the transition from the fast flow to the slow one happens at the points ${F=\textstyle\left(-4/\sqrt{3},\,16/3\sqrt{3}\right)}$ and ${\textstyle D=\left(4/\sqrt{3},\,-16/3\sqrt{3}\right)}$.
     
    Let $\gamma_0$ be the periodic trajectory of system \eqref{fhn} with initial condition at $D$ (cf. Figure \ref{fig:caso5extau}). As we have seen, $\gamma_0$ follows periodically the cycle:
    \begin{equation}
        \label{ciclo}
        \gamma_0: D\;\overset{\text{slow flow}}{\longrightarrow}\; P_+\; \overset{\text{fast flow}}{\longrightarrow}\; F\; \overset{\text{slow flow}}{\longrightarrow}\; P_- \;\overset{\text{fast flow}}{\longrightarrow}\; D\; \overset{\text{slow flow}}{\longrightarrow}\; \left(...\right)\;.
    \end{equation}

    We intend to compute the period $\mathcal{T}_{\gamma_0}$ of this cycle. Since $\varepsilon=0$, then the fast flow of the system is traversed almost instantaneously.
    Therefore, $\mathcal{T}_{\gamma_0}$ is determined by the durations $\mathcal{T}_{\gamma_0}^{\mathcal{C}_{0\text{L}}}$ and $\mathcal{T}_{\gamma_0}^{\mathcal{C}_{0\text{R}}}$ of $\gamma_0$ along the branches $\mathcal{C}_{0\text{L}}$ and $\mathcal{C}_{0\text{R}}$, respectively. Since $f$ is odd, we have $\mathcal{T}_{\gamma_0}^{\mathcal{C}_{0\text{L}}}=\mathcal{T}_{\gamma_0}^{\mathcal{C}_{0\text{R}}}$. Thus,
    \begin{equation}{\label{period_eps0}}
        \mathcal{T}_{\gamma_0}\approx \mathcal{T}_{\gamma_0}^{\mathcal{C}_{0\text{L}}} + \mathcal{T}_{\gamma_0}^{\mathcal{C}_{0\text{R}}}=2\;\mathcal{T}_{\gamma_0}^{\mathcal{C}_{0\text{R}}}.   
    \end{equation}
    Hence, the time spent by the trajectory in $\mathcal{C}_{0\text{R}}$ is given by: 
    \begin{equation*}
        \mathcal{T}_{\gamma_0}^{\mathcal{C}_{0\text{R}}} \equiv \int_0^{\mathcal{T}_{\gamma_0}^{\mathcal{C}_{0\text{R}}}} \; \mathrm{d}\tau.    
    \end{equation*}
    Using \eqref{expVC} we determine $\mathcal{T}_{\mathcal{C}_{0\text{R}}}$ with respect to the fast variable:
    \begin{equation}
        \label{periodo}
        \mathcal{T}_{\gamma_0} \approx 2\;\mathcal{T}_{\gamma_0}^{\mathcal{C}_{0\text{R}}} = 2 \int_{\frac{4}{\sqrt{3}}}^{\frac{2}{\sqrt{3}}} \dfrac{4-3x^2}{bx^3 + (1-4b)x + c}\; \mathrm{d}x ,
    \end{equation}
    that may be evaluated for any suitable value of $b$ and $c$. For instance, for $b=c=0$ we are in case (v) and
    \begin{equation*}
         \mathcal{T}_{\gamma_0} = 2\int_{\frac{4}{\sqrt{3}}}^{\frac{2}{\sqrt{3}}} \frac{4-3x^2}{x}\; \mathrm{d}x = 12 - 8\log\left(2\right) \approx 6.45\,.
    \end{equation*}
    %
    %
    %
    %
    %

    In this section we have described the dynamics of \eqref{fhn} in the singular case $\varepsilon=0$, in the various cases according to the number and stability of its equilibria, as shown in Figures \ref{fig:casos_1-4} and \ref{fig:caso5extau}. We have also obtained conditions under which \eqref{fhn} with $\varepsilon=0$ has a periodic solution that alternates between the fast and the slow flows and we have determined its period. In the next section we look at the dynamics of \eqref{deflento} in the nonsingular case, i.e., when $\varepsilon \gtrsim 0$.
    %
\section{The Regular Case}\label{sec:reg}
    %
    The independence  of the time scales when $\varepsilon=0$ makes the system simpler to understand.
    By allowing the dynamics to be separated into two distinct phases, slow and fast, we can study each one independently and understand its underlying dynamics.

    As we will see in Fenichel's theorem (Theorem \ref{thfenichel} stated below) when $0<\varepsilon\ll1$, the dynamics of the perturbed system is similar to that of the singular system, with a deviation of $\mathcal{O}\left(\varepsilon\right)$, where $\mathcal{O}$ stands for the usual \emph{Landau} notation.
    That is, the smaller the $\varepsilon$, the \emph{more similar} the system trajectories are to those described in the singular system  (cf. Figure \ref{fig:fhn11_regular_xy}). Fenichel's theorem guarantees that for $0<\varepsilon\ll1$, there exists a set $S_\varepsilon$ very close to $S_0$ that exhibits the {\emph same} behaviour as $S_0$.
     
  
    \begin{theorem}[Fenichel, 1979 \citeyear{FEN}]\label{thfenichel}
        Suppose $S_0$ is a compact normally hyperbolic submanifold of the critical manifold $\VC$ of \eqref{deflento} and that $f,g\in{C^r}\,\,\left(2\le r<\infty\right)$.
         Then for $\varepsilon>0$, sufficiently small, the following hold:
        \begin{enumerate}[\hspace{0.5cm} (H1)]
            \item There exists a locally invariant manifold $S_\varepsilon$ diffeomorphic to $S_0$.
            
            \item\label{item:H2}
            The set $S_\varepsilon$ has Hausdorff\footnote{The \emph{Hausdorff distance} between two nonempty sets $V,W\subset\RR^{m+n}$ is defined by $$d_H\left(V,W\right):=\max\left\{\sup_{v\in V}\inf_{w\in W} || v-w||,\,\sup_{w\in W}\inf_{v\in V} || v-w|| \right\}.$$} distance $\mathcal{O}\left(\varepsilon\right)$, as $\varepsilon \to 0$, from $S_0$.
            
            \item The flow in $S_\varepsilon$ converges to the slow flow, as $\varepsilon\to 0$.
            
            \item The set $S_\varepsilon$ is $C^r$--smooth.
            
            \item The set $S_\varepsilon$ is normally hyperbolic and has the same stability properties with respect to the fast variables as $S_0$.
            
            \item The set $S_\varepsilon$ is usually not unique. In regions that remain at a fixed distance from the topological boundary $\partial S_\varepsilon$ of $S_\varepsilon$, all manifolds satisfying (H1)--(H5) lie at a Hausdorff distance $\mathcal{O}\left(e^{-K/\varepsilon}\right)$ from each other for some  $K>0,$ where $K=\mathcal{O}\left(1\right)$.
        \end{enumerate}
    \end{theorem}

    The proof of this theorem is extensive and covers topics in perturbation theory that are beyond the scope of this work. 
    The  reader may find it in \citet{kuehn} or in its original version in \citet{FEN}.
    \begin{definition}\label{def:variedadalenta}
        The set $S_\varepsilon$ mentioned in Theorem \ref{thfenichel} is referred to as the \textbf{slow manifold} of \eqref{deflento}.
    \end{definition}
    %
\subsection*{Discussion of Fenichel's theorem}
    Theorem \ref{thfenichel} is central to the study of systems with multiple time scales.  Let us examine the role of each statement separately.
         
    Statement {\sl (H6)} is useful to simplify the language. In fact, there is not just one slow manifold $S_\varepsilon$.
    However, since all possible slow manifolds $S_\varepsilon^j,\, j\in\mathbb{N},$ are at a distance $\mathcal{O}\left(e^{-K/\varepsilon}\right)\overset{\varepsilon\to 0}{\longrightarrow}0$ from each other, we can simplify the discussion and consider only one slow manifold $S_\varepsilon$.
         
    The existence of a one-to-one correspondence between $S_\varepsilon$ and $S_0\subset\VC$, {\sl (H1)}, is important   for the other statements in the theorem.
    The set $S_\varepsilon$, being locally invariant, means that trajectories do not leave the slow manifold except at its boundary. 
    Statement {\sl (H2)} guarantees that the smaller $\varepsilon$ is, the closer $S_\varepsilon$ is to $S_0$. Thus, since $S_\varepsilon$ is locally invariant, the trajectories on $S_\varepsilon$ approach $S_0$ by {\sl (H3)}.
    Finally, the fact that $S_\varepsilon$ is continuously differentiable to the same order as $S_0$, {\sl (H4)}, and is diffeomorphic to $S_0$, {\sl (H1)}, allows us to assert that the two sets $S_0$ and $S_\varepsilon$ share the same stability.    The set $S_\varepsilon$ being normally hyperbolic, {\sl (H5)}, ensures that the fast variables have the same stability properties of $S_0$.
    %
\subsection*{The manifold $S_\varepsilon$}
\label{sec:se}
    In this section, we give a precise way to describe $S_\varepsilon$. Consider the (1,1)-fast--slow system in the slow time scale:
    \begin{equation}\label{eq:slowTime}
        \begin{aligned}
            &\varepsilon\,\dfrac{\mathrm{d}x}{\mathrm{d}\tau}=\varepsilon\Dot{x}=f(x,y,\varepsilon),\\
            &\phantom{\varepsilon\,}\dfrac{\mathrm{d}y}{\mathrm{d}\tau}=\phantom{\varepsilon}\Dot{y}=g(x,y,\varepsilon),
        \end{aligned}
    \end{equation}         \noindent with $0<\varepsilon\ll1$.

    Let $S_0$ be a compact, normally attracting subset of $\VC$ without equilibria ($g\left(x,y,\varepsilon\right)\neq0$). By the implicit function theorem, there exists a function $x=x_0\left(y\right)$, for $y\in\left]a,b\right[$, such that we define $S_0$ as:
    \begin{equation}
        \begin{aligned}
            S_0&=\left\{\left(x,y\right)\in{\RR^2}:f\left(x,y,0\right)=0 \wedge y\in\left]a,b\right[\right\}\\
            &=\left\{\left(x,y\right)\in{\RR^2}:x=x_0\left(y\right) \wedge y\in\left]a,b\right[\right\}.
        \end{aligned}
    \end{equation}

    \begin{theorem}[Fenichel, 1979 \citeyear{FEN}. See also Theorem 11.1.1 in Kuehn, 2015 \citeyear{kuehn}.]
        Let $S_0$ be a compact normally hyperbolic submanifold of $\VC$. Then there exists a slow manifold $S_\varepsilon$ that is $\mathcal{O}\left(\varepsilon\right)$-close to $S_0$ for $\varepsilon>0$ sufficiently small. Locally, $S_\varepsilon$ is represented as a graph of a smooth function $h(y,\varepsilon)$:
        \begin{equation*}
            S_\varepsilon=\left\{\left(x,y\right)\in\RR^2:x= h(y,\varepsilon)\right\},
        \end{equation*}              %
        where the map $ y\mapsto h(y,\varepsilon):\RR\to\RR$  has the asymptotic expansion 
        \begin{equation*}
            h(y,\varepsilon)
            =h_0\left(y\right) + h_1\left(y\right)\varepsilon + h_2\left(y\right)\varepsilon^2 + {\mathcal O}(\varepsilon^3).
        \end{equation*}
    \end{theorem}
    \bigbreak
    From the previous result, we may conclude that $h_0(y)=x_0(y)$. The curve $h_0\left(y\right)$  defines locally the set $S_0$ of the critical manifold $\VC$.\\
        
    In what follows, our aim is to show how to compute the term $h_1(y)$ in the expansion. Taking into account that $g\left(x,y,\varepsilon\right)\neq0$, we write \eqref{eq:slowTime} as
    \begin{equation}\label{eq:dxdy}
        \varepsilon\,\dfrac{\mathrm{d}x}{\mathrm{d}y}=\dfrac{f\left(x,y,\varepsilon\right)}{g\left(x,y,\varepsilon\right)}. 
    \end{equation}
   
    By invariance, the set $x\equiv h(y,\varepsilon)$ satisfies \eqref{eq:dxdy} and hence
    \begin{equation}\label{eq:dhdy}
        \varepsilon \dfrac{dh_0}{dy}(y)+
        \varepsilon^2 \dfrac{dh_1}{dy}(y)+ {\mathcal O}(\varepsilon^3)=
        \dfrac{f\left(h(y,\varepsilon),y,\varepsilon\right)}{g\left(h(y,\varepsilon),y,\varepsilon\right)} .
    \end{equation}

    We want to estimate the first terms in \eqref{eq:dhdy}.
    For the first one, if $\left(x\left(\tau\right),y\left(\tau\right)\right)$ is a solution of \eqref{eq:slowTime} for $\varepsilon=0$ with initial condition in $S_0$, hence
    $f\left(x(\tau),y(\tau),0\right)\equiv 0$, then
    \begin{equation*}
        f_x\left(x,y,0\right) \dfrac{d x}{d \tau} + f_y\left(x,y,0\right) \dfrac{d y}{d \tau}=0
        \quad\Rightarrow\quad
        \dfrac{dx}{dy}=-\dfrac{f_y\left(x,y,0\right)}{f_x\left(x,y,0\right)}
    \end{equation*}
    which implies (because $x\equiv h(y,\varepsilon)$)
    \begin{equation}\label{eq:dh0}
        \dfrac{\partial h}{\partial y}(y,0)=\dfrac{d h_0}{d y}(y)=
        -\dfrac{f_y\left(h_0(y),y,0\right)}{f_x\left(h_0(y),y,0\right)}.
    \end{equation}

    Now we compute $h_1(y)$. 
    To do this, let $H(x,y,\varepsilon):=\dfrac{f(x,y,\varepsilon)}{g(x,y,\varepsilon)}$, so the right hand side of \eqref{eq:dhdy} is 

    \begin{equation}\label{eq:Hepsilon2}
        H\left(h(y,\varepsilon),y,\varepsilon\right)=
        \varepsilon\left(\dfrac{\partial H}{\partial x}\left(h(y,0),y,0\right) h_1(y)+ 
        \dfrac{\partial H}{\partial \varepsilon}\left(h(y,0),y,0\right)  \right) + {\mathcal O}(\varepsilon^2) ,
    \end{equation}
    since $H\left(h(y,0),y,0\right)=0$ and $\dpt \frac{\partial h}{\partial \varepsilon}(y,0)=h_1(y)$.\\
    
    Next, we use again $f(h(y,0), y, 0)=0$ to  compute the derivatives of $H$ appearing in \eqref{eq:Hepsilon2}:


    \begin{equation*} \label{eq:dH}
        \dfrac{\partial H}{\partial x}\left(h(y,0),y,0\right) = \dfrac{f_x\left(h(y,0),y,0\right) }{g\left(h(y,0),y,0\right)}, \qquad 
        \dfrac{\partial H}{\partial\varepsilon}\left(h(y,0),y,0\right) = \dfrac{f_\varepsilon\left(h(y,0),y,0\right)}{g\left(h(y,0),y,0\right)} \\
    \end{equation*}       
    and \eqref{eq:Hepsilon2} becomes
    \begin{equation}\label{eq:Hepsilon3}
        H\left(h(y,\varepsilon),y,\varepsilon\right)=
        \varepsilon\left(\dfrac{f_x\left(h_0(y),y,0\right) }{g\left(h_0(y),y,0\right)} h_1(y)+ 
        \dfrac{f_\varepsilon\left(h_0(y),y,0\right)}{g\left(h_0(y),y,0\right) } \right) + {\mathcal O}(\varepsilon^2) .
    \end{equation}

    Substituting \eqref{eq:dh0} and \eqref{eq:Hepsilon3} into \eqref{eq:dhdy} yields
    \begin{equation*}\label{eq:dhdy2}
        -\varepsilon\dfrac{f_y\left(h_0(y),y,0\right)}{f_x\left(h_0(y),y,0\right)}+ {\mathcal O}(\varepsilon^2) =
        \varepsilon\left(\dfrac{f_x\left(h_0(y),y,0\right) }{g\left(h_0(y),y,0\right)} h_1(y)+ 
        \dfrac{f_\varepsilon\left(h_0(y),y,0\right)}{g\left(h_0(y),y,0\right) } \right) + {\mathcal O}(\varepsilon^2)
    \end{equation*}
    hence 
    \begin{equation}\label{eq:2h1}
        h_1(y)= -\dfrac{f_y\left(h_0(y),y,0\right)g(h_0(y),y,0)}{\left(f_x\left(h_0(y),y,0\right)\right)^2}-\dfrac{f_\varepsilon\left(h_0(y),y,0\right)}{f_x\left(h_0(y),y,0\right) }.
    \end{equation}

    \bigbreak

\begin{figure}
        \centering
        \begin{subfigure}{0.3\linewidth}
            \centering
            \includegraphics[width=\linewidth]{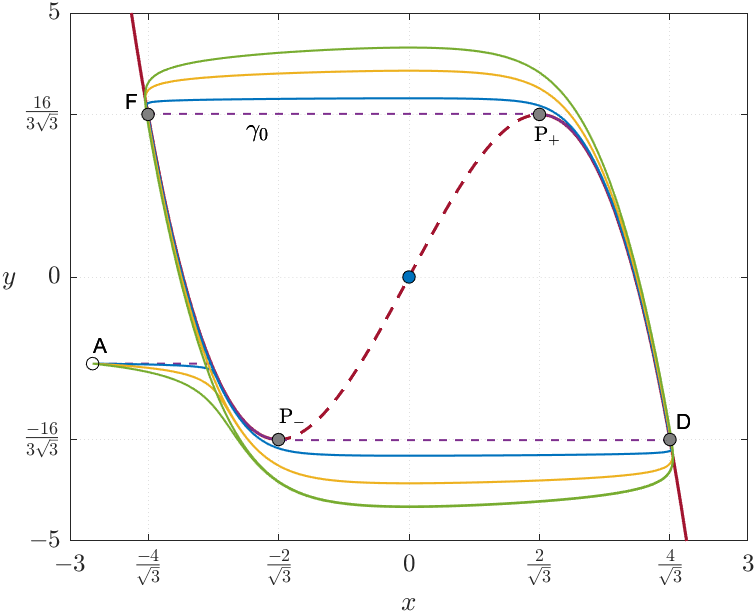}
            \caption{ \small }
            \label{fig:fhn11_regular_xy}
        \end{subfigure}
        \begin{subfigure}{0.3\linewidth}
            \centering
            \includegraphics[width=\linewidth]{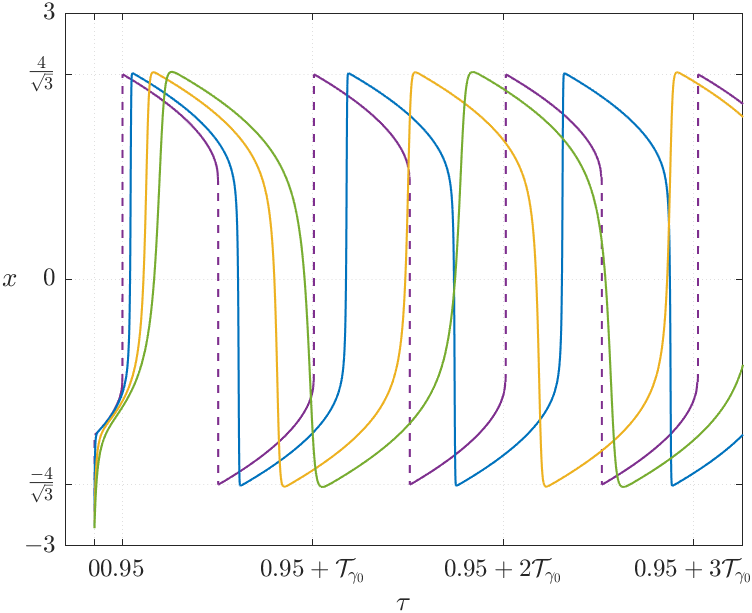}
            \caption{ \small }
            \label{fig:fhn11_regular_xtau}
        \end{subfigure}
        \begin{subfigure}{0.3\linewidth}
            \centering
            \includegraphics[width=\linewidth]{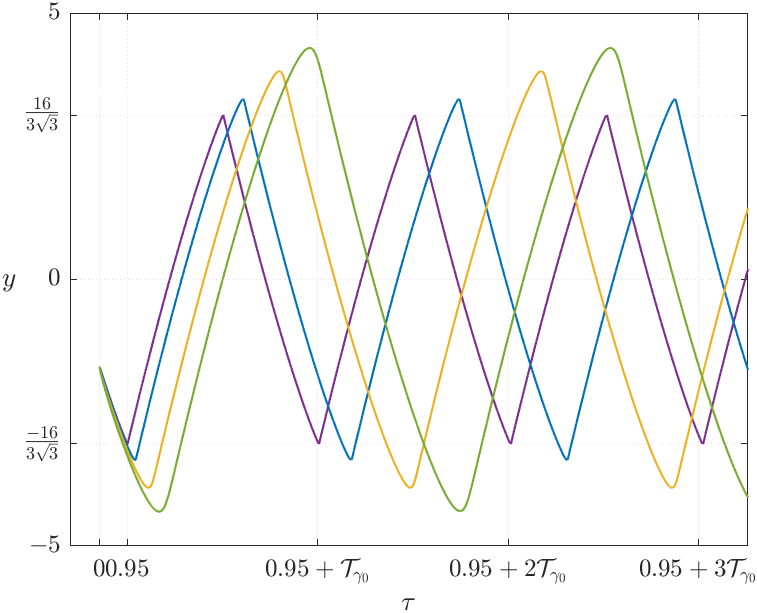}
            \caption{ \small }
            \label{fig:fhn11_regular_ytau}
        \end{subfigure}
        \caption{ \small  (a) - Trajectory of the system \eqref{fhn} with parameters $b=0=c$, started at $A\approx\left(-2.8,1.64\right)$ and with $\varepsilon\in\left\{0,0.1,0.5,1\right\}$ (violet, blue, green, yellow). Critical manifold (red) and unstable equilibrium point (blue). (b) - Dynamics of the trajectories of (a) in the coordinate plane $\left(\tau,x\right)$. (c) - Dynamics of the trajectories of (a) in the coordinate plane $\left(\tau,y\right)$.}
        \label{fig:fhn11_regular}
    \end{figure}

    \begin{example}
        Consider the system FH-N \eqref{fhn}. If we take $S_0\subset\mathcal{C}_{0\text{L}}$, we get that $S_0$ is defined as 
        \begin{equation*}
            h_0\left(y\right)=-\dfrac{4\left(\dfrac{2}{3}\right)^{1/3}}{\left(9y+\sqrt{3\left(27y^2-256\right)}\right)^{1/3}}-\dfrac{\left(9y+\sqrt{3\left(27y^2-256\right)}\right)^{1/3}}{18^{1/3}}.     
        \end{equation*}
        Additionally, 
        using \eqref{eq:2h1}, we get
        \begin{equation*}
            h_1\left(y\right)=\dfrac{h_0\left(y\right)-by-c}{\left(4-h_0^2\left(y\right)\right)^2} 
        \end{equation*}
        and $S_\varepsilon$ is approximated by the graph of
        \begin{equation*}
            h_\varepsilon\left(y\right)=h_0\left(y\right) \,+\, \dfrac{h_0\left(y\right)-by-c}{\left(4-h_0^2\left(y\right)\right)^2} \,\varepsilon. \qquad \diamond    
        \end{equation*}
    \end{example}
    Since we defined the slow manifold $S_\varepsilon$ as an asymptotic expansion starting at the critical manifold $\VC$, we can apply the same reasoning to determine the period of a periodic solution, $\mathcal{T}_{\gamma_\varepsilon}$, of a fast-slow system, with $0<\varepsilon\ll1$  (cf. Figure \ref{fig:fhn11_regular_xtau}).
    In particular, using the same argument used to prove \eqref{period_eps0} we have:
    \begin{equation*}
        \mathcal{T}_{\gamma_\varepsilon} = \mathcal{T}_{\gamma_0} + \mathcal{O}\left(\varepsilon\right).    
    \end{equation*}  
     
    %
\section{ Bifurcations}\label{sec:bif}
    To understand the effect of the parameters $b$ and $c$, we analyse two cases of \eqref{fhn}, where the first has $b=0$ and $c\in\RR$, and the second has $b\in\RR\setminus{\left\{0\right\}}$ and $c=0$. \\
    
    (i) For $c\in\RR$ and $b=0$, system \eqref{fhn} may be recast into the form:
        \begin{equation}
            \label{fhn_b=0}
            \begin{aligned}
                &\dfrac{\mathrm{d}x}{\mathrm{d}t}=x'=-y+4x-x^3\\
                &\dfrac{\mathrm{d}y}{\mathrm{d}t}=y'=\varepsilon (x-c),
            \end{aligned}
        \end{equation}
        with $0<\varepsilon\ll1$. It has only one equilibrium $\mathrm{E}=\left(c,\varphi(c)\right)$ for $\varphi(x)\coloneqq4x-x^3$, which results from the intersection of the cubic curve $y=\varphi(x)$ and the  {vertical} line $x=c$.\\
        
        In order to evaluate the stability of $\mathrm{E}$, we compute the jacobian matrix $J_E$ of the vector field associated to \eqref{fhn_b=0} at $E$ and find its eigenvalues, say $\lambda_\mathrm{E}^\pm\left(c,\varepsilon\right)$. Thus, we get
        \begin{align*}
            J_{\mathrm{E}}= \begin{bmatrix}
                                      \, 4-3c^2 \, \phantom{aaa} \, & -1\\
                                    \varepsilon & 0
                            \end{bmatrix} \qquad \text{and} \qquad \lambda_\mathrm{E}^\pm\left(c,\varepsilon\right) = \dfrac{\left(4-3c^2\right)\pm\sqrt{\left(4-3c^2\right)^2-4\varepsilon}}{2}.
        \end{align*}
        It is straightforward to conclude that:
        
        \begin{enumerate}[\hspace{0.5cm} (1)]
            \item The equilibrium $\mathrm{E}$ is an unstable node/focus if $|c|<\dfrac{2}{\sqrt{3}}$ and 
            \item $\mathrm{E}$ is a stable node/focus if $|c|>\dfrac{2}{\sqrt{3}}$.
        \end{enumerate}
         
\begin{figure}[h!]
            \centering
            \begin{subfigure}{0.49\linewidth}
                \centering
                \includegraphics[width=\linewidth]{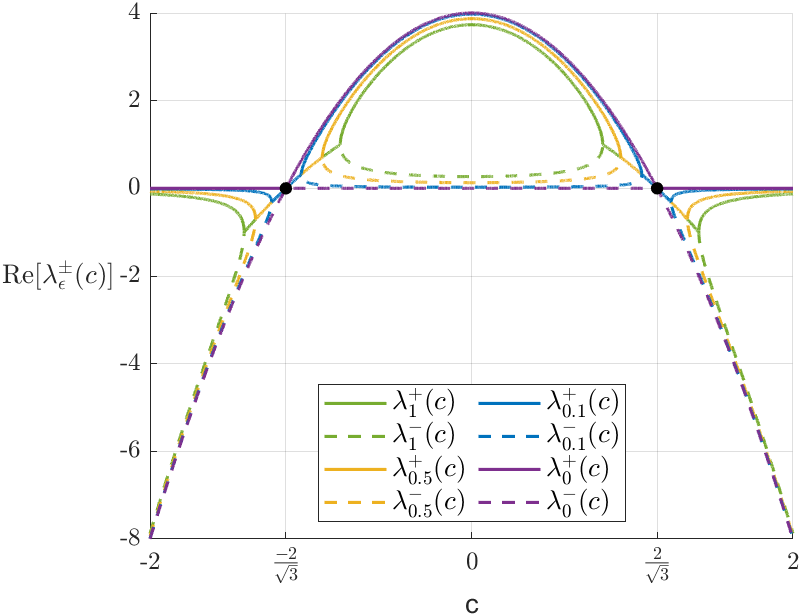}
                \caption{ \small  $\mathrm{Re}\left[\lambda_\varepsilon^\pm\left(c\right)\right]$}
                \label{fig:RebifurC}
            \end{subfigure}
            \begin{subfigure}{0.49\textwidth}
                \centering
                \includegraphics[width=\linewidth]{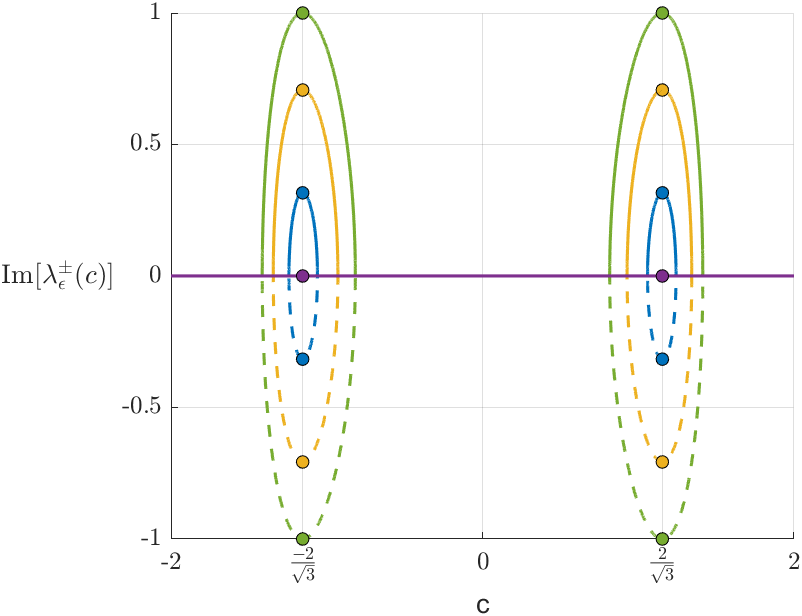}
                \caption{ \small  $\mathrm{Im}\left[\lambda_\varepsilon^\pm\left(c\right)\right]$}
                \label{fig:ImbifurC}
            \end{subfigure}
            \caption{ \small Evolution of the real and imaginary parts of the eigenvalues of the Jacobian matrix of system \eqref{fhn_b=0}, evaluated at the equilibrium point $E$, as a function of the parameter $c$ ($\lambda_\varepsilon^+$ with solid line and $\lambda_\varepsilon^-$ with dashed line) for different values of $\varepsilon$, $\varepsilon=1$ (green), $\varepsilon=0.5$ (yellow), $\varepsilon=0.1$ (blue) e $\varepsilon=0$ (violet). Hopf bifurcation (black).}
            \label{fig:bifurC_lambdas}
        \end{figure}

        When $c=\pm\mathrm{c^H}\text{, with } \mathrm{c^H} = 2/\sqrt{3}$, and $0<\varepsilon\ll1$, we get $\lambda_\mathrm{E}^\pm\left(\mathrm{c},\varepsilon\right) = \pm\mathrm{i}\sqrt{\varepsilon}$.
        For $\varepsilon>0$, as shown in Figure \ref{fig:bifurC_lambdas}, the stability of the equilibrium point in the neighbourhood of $\mathrm{c^H}$ changes from a stable focus to an unstable focus as $c$ decreases.    Similarly, by symmetry, the dual stability transition occurs in the neighbourhood of $-\mathrm{c^H}$ as $c$ increases.
         
\begin{figure}
            \centering
            \begin{subfigure}{0.46\linewidth}
                \includegraphics[width=\linewidth]{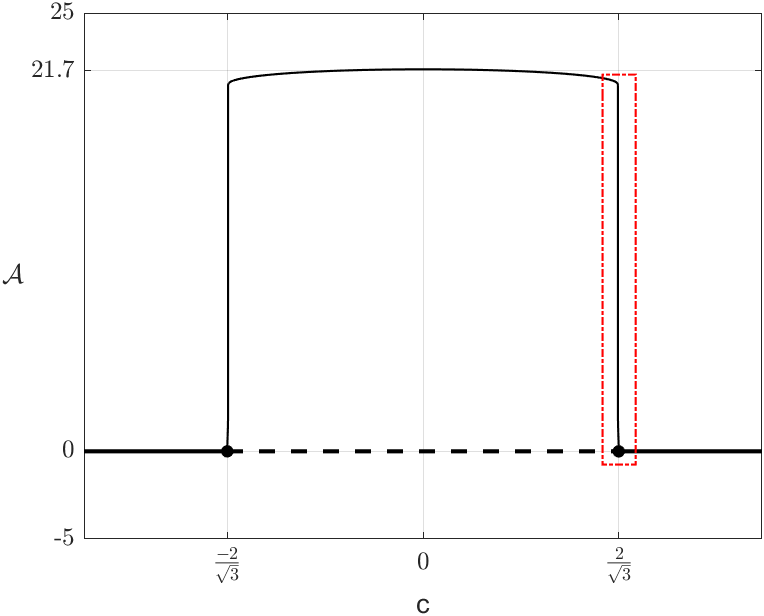}
                \caption{ \small }
                \label{fig:bifc_eps05_diagram}
            \end{subfigure}
            \hfill
            \begin{subfigure}{0.46\linewidth}
                \includegraphics[width=\linewidth]{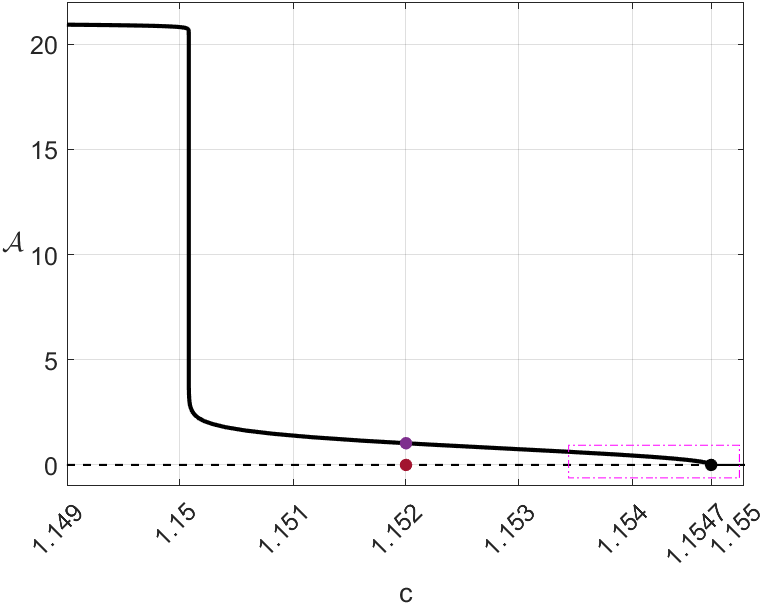}
                \caption{ \small }
                \label{fig:bifc_eps05_diagram_amp}
            \end{subfigure}
            \vfill
            \begin{subfigure}{0.46\linewidth}
                \includegraphics[width=\linewidth]{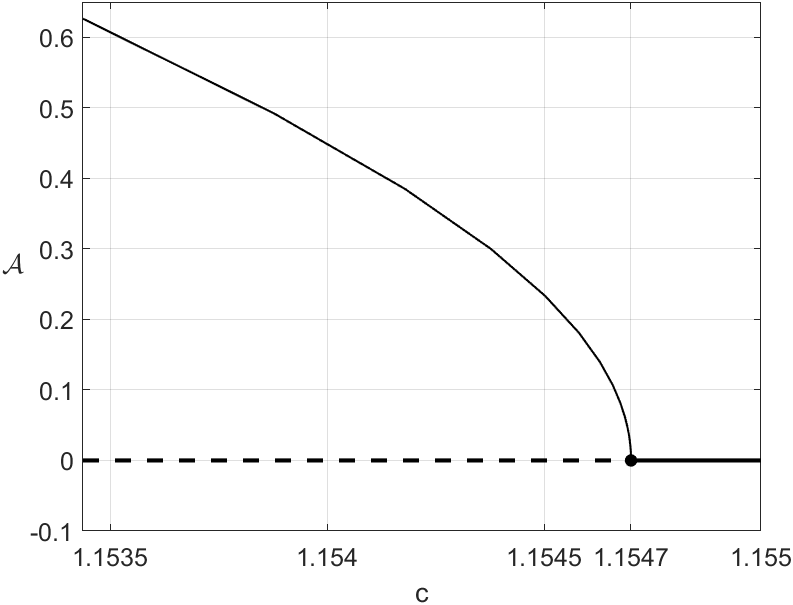}
                \caption{ \small }
                \label{fig:bifc_eps05_diagram_2xamp}
            \end{subfigure}
            \hfill
            \begin{subfigure}{0.46\linewidth}
                \includegraphics[width=\linewidth]{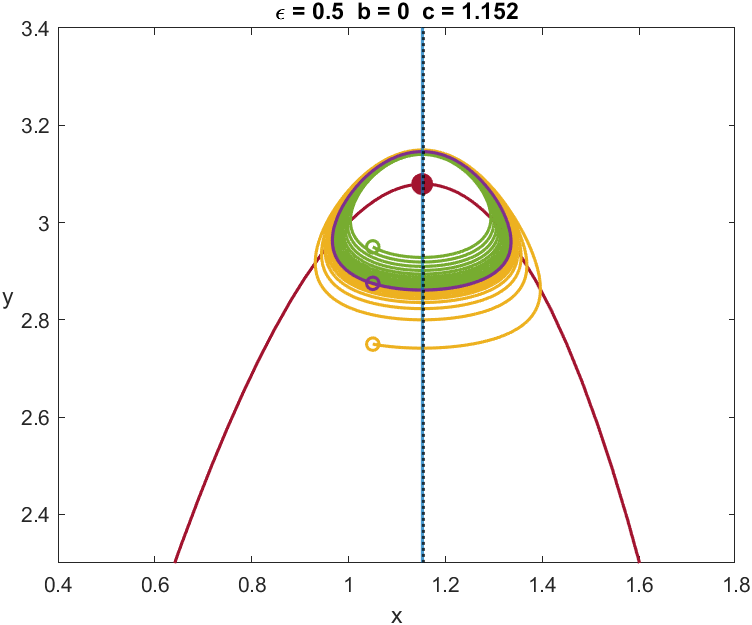}
                \caption{ \small }
                \label{fig:b0_c1152_e05}
            \end{subfigure}
            \caption{ \small  (a), (b) and (c) Bifurcation diagrams: length ${\mathcal A}$ of the periodic solutions of system \eqref{fhn_b=0} with $\varepsilon=0.5$ as a function of the parameter $c$ with different levels of magnification. Equilibria correspond to ${\mathcal A}=0$. The graphs in (b) and (c) are zooms of the red and magenta rectangles in (a) and (b), respectively. Conventions: Stable periodic solutions (solid line at $\mathcal{A}\neq0$) for $|c|<2/\sqrt{3}$. (d) Example of a periodic solution of the system \eqref{fhn_b=0} for $c=1.152$ and $\varepsilon=0.5$. $\VC$ (red), vertical lines $x=c$ (blue) and $x=2/\sqrt{3}$ (black). Approaching/leaving trajectories from the equilibrium point (yellow/green) and periodic solution (violet). }
            \label{fig:bifc_eps05_mais_exemplo}
        \end{figure}

        For $c=1.15<\mathrm{c^H}$ and $\varepsilon=1$, the equilibrium point $\mathrm{E}$ is an unstable focus, and thus it is expected that trajectories starting close to the fold point $P_+$ diverge, for $t>0$.
        Numerically, we found a periodic solution when $c\lessapprox \mathrm{c^H}$ (cf. Figure \ref{fig:b0_c1152_e05}).
        This dynamics is typical of a \emph{subcritical Hopf bifurcation}, with a stable periodic solution  as described in the next result. Following the terminology of \cite[Cap. IV, Section 2]{schaeffer1985singularities}, \emph{subcritical} means that the periodic solution exists for $c<c^\mathrm{H}$. If the periodic solutions are found for $c>c^\mathrm{H}$, the bifurcation would be called \emph{supercritical}.

        \begin{theorem}[Hopf, 1942 \cite{Hopf}. See also Marsden \& McCraken, 1976 \citeyear{marsden} for English translation in Section 5]\label{th:Hopf}
            Given a one-parameter family of differential equations in $\RR^2$ with parameter $c$:
            \begin{align*}
                \Dot{x}&=f(x,y;c)\\
                \Dot{y}&=g(x,y;c), \qquad f,g\in{{C^k}}\quad \mathrm{(k}\geq 2\mathrm{)},
            \end{align*}
            such that there exists an equilibrium point $\left(x_0\left(c\right),y_0\left(c\right)\right)\in{\RR^2}$ for all $c\in\RR$ and the eigenvalues of the jacobian matrix $J_{\left(x_0\left(c\right),y_0\left(c\right)\right)}$ can be written as $\lambda\left(c\right)=\alpha\left(c\right)  \pm    i\,   \beta \left(c\right)$,
            if for some $c^*\in\RR$ the following holds: 
            \begin{equation*}
                \alpha\left(c^*\right)=0, \quad \beta\left(c^*\right)\neq0 \quad \text{ and } \quad{\dfrac{\mathrm{d}\alpha}{\mathrm{d}c}}({c^*})\neq 0,    
            \end{equation*}
            then, for $c$ close to $c^*$, there exists \textbf{a periodic solution} (limit cycle) around $\left(x_0\left(c\right),y_0\left(c\right)\right)$.
        \end{theorem}
        %
        {System \eqref{fhn_b=0}  satisfies the conditions of Theorem \ref{th:Hopf},
        thus confirming the numerical finding of Figure \ref{fig:b0_c1152_e05}.   The direction of bifurcation and the stability of the periodic solution will be discussed in more detail in Section \ref{sec:canards} below.}

\begin{figure}
            \centering
            \begin{subfigure}{0.32\linewidth}
                \includegraphics[width=\linewidth]{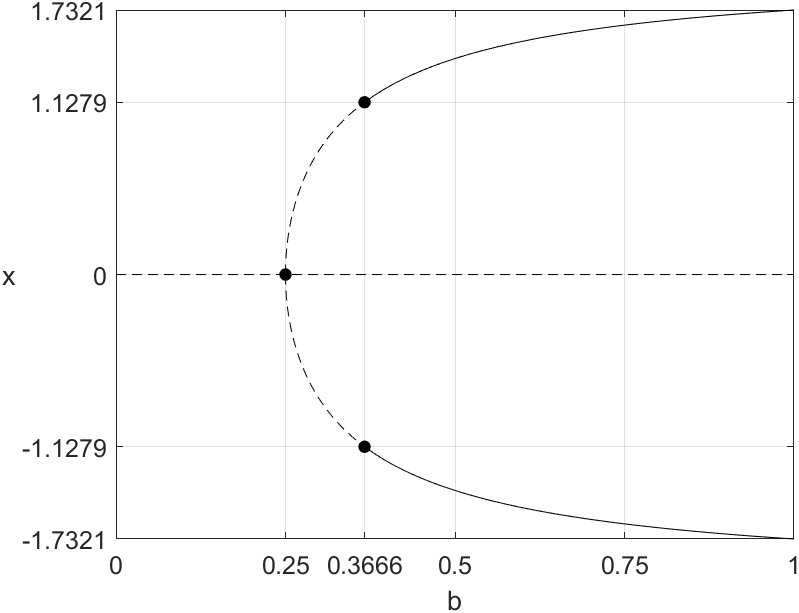}
                \caption{ \small }
                \label{fig:bifb_eps05_diagram}
            \end{subfigure}
            \begin{subfigure}{0.32\linewidth}
                \includegraphics[width=\linewidth]{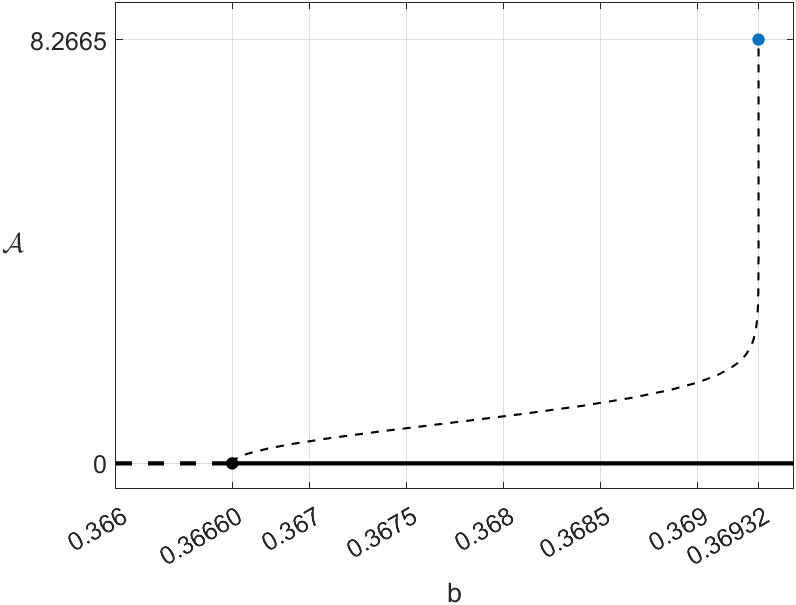}
                \caption{ \small }
                \label{fig:bifb_eps05_diagram_Eq+_amp}
            \end{subfigure}
            \begin{subfigure}{0.34\linewidth}
                \includegraphics[width=\linewidth]{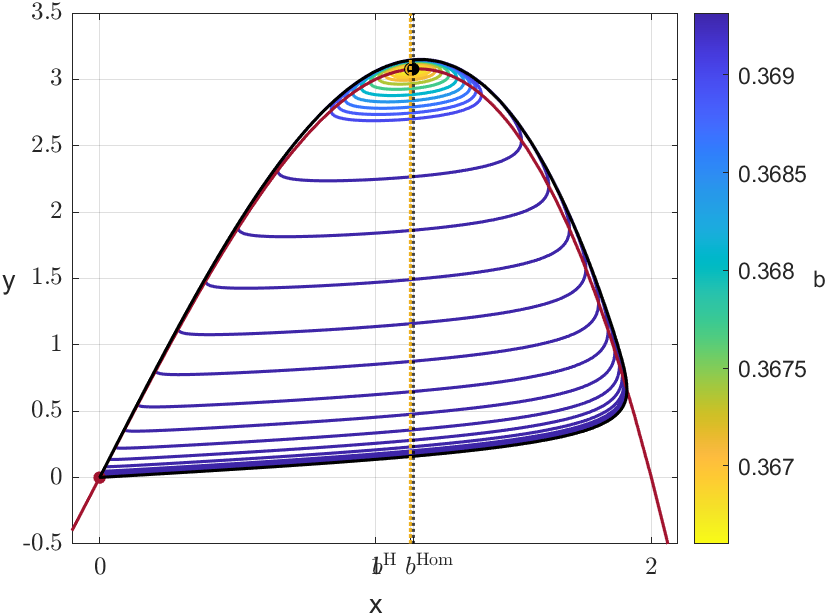}
                \caption{ \small }
                \label{fig:bifb_eps05_InstOrb_HomOrb}
            \end{subfigure}
            \caption{ \small  (a) - Bifurcation of the equilibria of the system \eqref{fhn_c=0} as a function of the parameter $b$, with $\varepsilon=0.5$. Equilibrium point $\mathrm{E}_0$ is always unstable (dashed line). Equilibria $\mathrm{E}_-$ and $\mathrm{E}_+$ are unstable for $1/4<b<b_{0.5}^H\approx 0.3666$ and stable $b>b_{0.5}^H$ (solid line). Bifurcation points at $b=1/4$ (Pitchfork) and at $b=b_{0.5}^H$ (supercritical Hopf). (b) Length of the periodic solution around the point $\mathrm{E}_+$ as a function of $b$. Supercritical Hopf bifurcation (black dot) leading to an unstable periodic solution that collides with the homoclinic orbit associated to the saddle point $\mathrm{E}_0$. Homoclinic bifurcation at $b_{0.5}^{Hom}\approx 0.36932$ (blue dot). (c) - Example of unstable periodic solutions of the system \eqref{fhn_c=0} generated by the Hopf bifurcation, for $b_{0.5}^H\leq b \leq b_{0.5}^{Hom}$ and $\varepsilon=0.5$. Homoclinic orbit (black line). Critical manifold $\VC$ and equilibrium $\mathrm{E}_0$ in red, vertical lines $x=x\left(b_{0.5}^H\right)$ (yellow) and $x=x\left(b_{0.5}^{Hom}\right)$ (black). Colour scale for the variation of the parameter $b$   between values $b_{0.5}^H$ and $b_{0.5}^{Hom}$}.
            \label{fig:bifb_eps05_mais_exemplo}
        \end{figure}

    (ii) For $b\in\RR\setminus{\left\{0\right\}}$ and $c=0$, system \eqref{fhn} may be written as:
        \begin{equation}
            \label{fhn_c=0}
            \begin{aligned}
                &\dfrac{\mathrm{d}x}{\mathrm{d}t}=x'=-y+4x-x^3\\
                &\dfrac{\mathrm{d}y}{\mathrm{d}t}=y'=\varepsilon (x-by),
            \end{aligned}
        \end{equation}
        with $0<\varepsilon\ll1$. For all $b\in\RR\setminus{\left\{0\right\}}$, the origin, $\mathrm{E}_0 =\left(0,0\right)$, is an equilibrium point the system \eqref{fhn_c=0}.
        If  either $b<0$ or $b>1/4$, the system has two symmetric equilibria at $x_\pm=\pm\sqrt{4-b^{-1}}$, as a result of a Pitchfork bifurcation, as the reader may check in Figure  \ref{fig:bifb_eps05_diagram}.
        They are: $$\mathrm{E}_\pm =\left(x_\pm,\varphi(x_\pm)\right), \quad \text{with}\quad  \varphi(x)\coloneqq 4x-x^3.$$
     
        For $b<0$, the equilibrium $E_0$ is an unstable node/focus and $E_-$ and $E_+$ are saddles.
        For $0<b\leq 1/4$, there is a stable limit cycle ({similar scenario to that of (v) represented in Figure \ref{fig:caso5}}) and only the origin is an equilibrium point of the system, but it still exhibits the same unstable node/focus dynamics.
     
\begin{figure}
            \centering
            \begin{subfigure}{.575\linewidth}
                \centering
                \includegraphics[width=\linewidth]{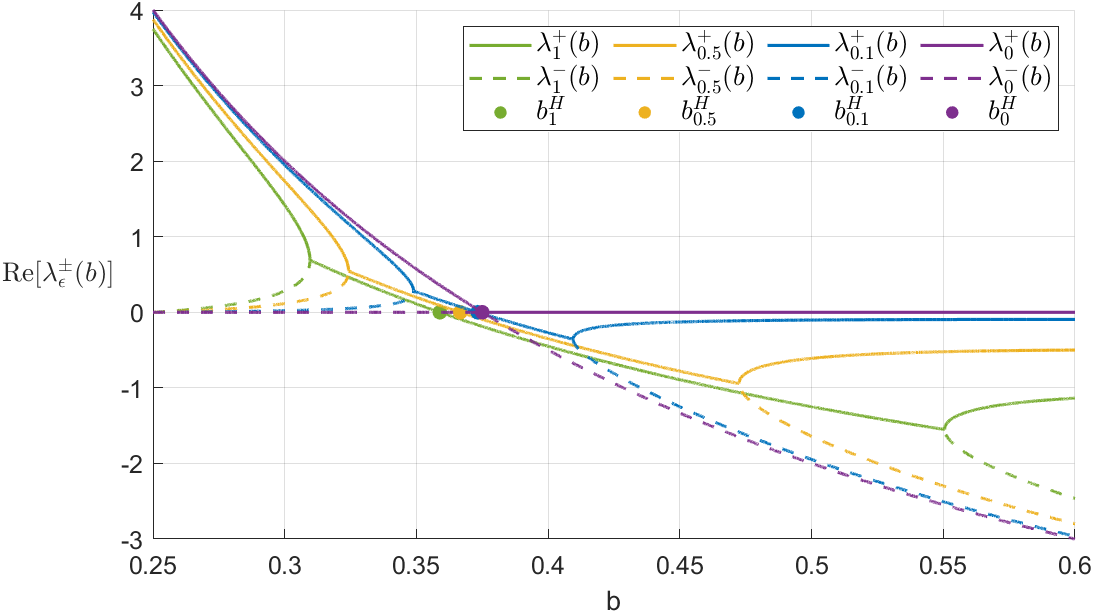}
                \caption{ \small  $\mathrm{Re}\left[\lambda_\varepsilon^\pm\left(b\right)\right]$.}
                \label{fig:RebifurB}
            \end{subfigure}%
            \begin{subfigure}{.425\textwidth}
                \centering
                \includegraphics[width=\linewidth]{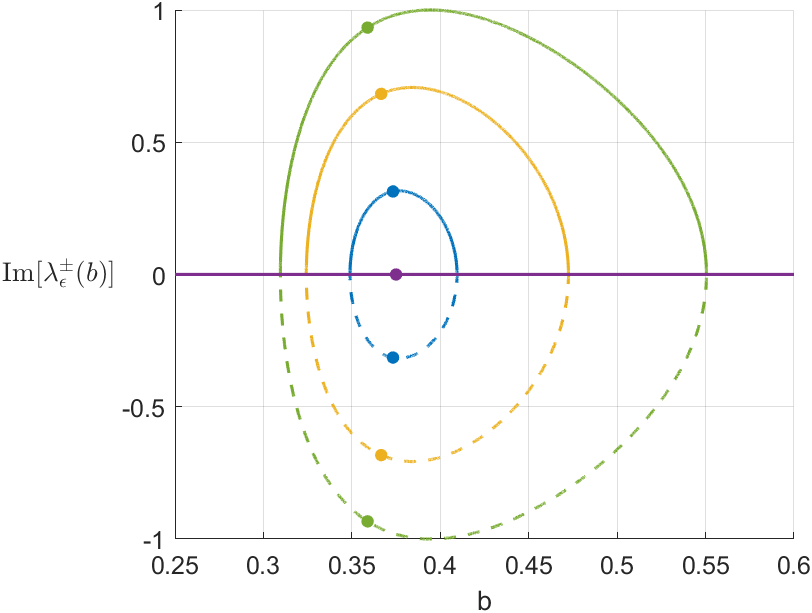}
                \caption{ \small  $\mathrm{Im}\left[\lambda_\varepsilon^\pm\left(b\right)\right]$.}
                \label{fig:ImbifurB}
            \end{subfigure}
            \caption{ \small  Evolution of the real and imaginary parts of the eigenvalues of the Jacobian matrix of the system \eqref{fhn_c=0}, evaluated at the equilibrium point, as function of the parameter $b$ ($\lambda_\varepsilon^+$ with solid line, $\lambda_\varepsilon^-$ with dashed line, and Hopf bifurcation point) for different values of $\varepsilon$, $\varepsilon=1$ (green), $\varepsilon=0.5$ (yellow), $\varepsilon=0.1$ (blue) e $\varepsilon=0$ (violet).}
            \label{fig:bifurB_lambdas}
        \end{figure}

        The equilibrium $\mathrm{E}_0$ becomes a saddle when $b>1/4$. Moreover, at $b=1/4$ two new equilibria $\mathrm{E}_-$ and $\mathrm{E}_+$  are created (cf. Figure \ref{fig:bifb_eps05_diagram}).  Given that $E_+$ and $E_-$ are symmetric, the stability of both equilibrium points is the same. Thus, we will only study the stability of $\mathrm{E}_+$ by computing  
        \begin{align*}
            J_{\mathrm{E}_+}= \begin{bmatrix}
                                    3\, b^{-1} -8 \phantom{aaa} & -1\\
                                    \varepsilon & -\varepsilon\, b
                                \end{bmatrix}, \quad \Tr\left(J_{\mathrm{E}_+}\right)=-\varepsilon\, b+3\, b^{-1} -8 \quad \text{and} \quad \det\left(J_{\mathrm{E}_+}\right)=2\varepsilon\left(4b-1\right),
        \end{align*}
        where $\Tr$ and $\det$ denote the trace and determinant operators. Note that $\det\left(J_{\mathrm{E}_+}\right)>0$, since $b>1/4$. The eigenvalues of $J_{\mathrm{E}_+}$ satisfy
        \begin{equation*}
            \lambda_{\varepsilon}^\pm\left(b\right) = \dfrac{\Tr\left(J_{\mathrm{E}_+}\right)\pm\sqrt{\Tr^{\,2}\left(J_{\mathrm{E}_+}\right)-4\,\det\left(J_{\mathrm{E}_+}\right)}}{2}.    
        \end{equation*}
        If we take
        \begin{equation*}
            \mathrm{b}_\varepsilon^\mathrm{H}=\dfrac{-4+\sqrt{16+3\varepsilon}}{\varepsilon}>\dfrac{1}{4},\qquad0<\varepsilon\ll1    
        \end{equation*}
        then at $b=\mathrm{b}_\varepsilon^\mathrm{H}$
        \begin{equation*}
            \Tr\left(J_{\mathrm{E}_+}\right) =0 \quad \text{and}\quad \eval{\dfrac{\mathrm{d}\mathrm{Re}\left(\lambda_{\varepsilon}^+\right)}{\mathrm{d}\,b}}_{b=\mathrm{b}_\varepsilon^\mathrm{H}}<0    
        \end{equation*}
        and therefore, by Theorem \ref{th:Hopf}, there is a Hopf bifurcation.
        In the limit, as $\varepsilon\to 0$, we have $\mathrm{b}_\varepsilon^\mathrm{H}\to 3/8$ (cf. Figure \ref{fig:bifurB_lambdas}).      
        
        Moreover, for $b=3/8$, the fold points are equilibria of \eqref{fhn_c=0}. \\
    
            The numerical simulations shown in Figure \ref{fig:bifb_eps05_mais_exemplo} indicate that the Hopf bifurcation at $\mathrm{E}_+$ is supercritical, hence generating an unstable periodic solution.   This will be confirmed analytically in Section \ref{sec:canards} below. This means that all trajectories starting in the regions bounded by the periodic solutions are attracted to the equilibrium point $\mathrm{E}_-$ or $\mathrm{E}_+$, also lying in the interior of that region. On the other hand, when trajectories start outside these regions, they converge to the global stable limit cycle which already existed when the equilibria were unstable. As observed in the example of Figure \ref{fig:bifb_eps05_mais_exemplo}, the  length  (a numerical estimate of the perimeter of the curve) of the unstable periodic solutions grows very quickly for small variations in $b$, but reaches an unexpected maximum for $b=b_{\varepsilon}^{Hom}$.   This happens when the periodic solution captures the saddle equilibrium $\mathrm{E}_0$ at a \emph{\textbf{homoclinic bifurcation}} \cite{guk,strogatz}.
        \begin{definition}
            Consider a system of  ordinary differential equations $\dfrac{\mathrm{d}X}{\mathrm{d}t}=F(X)$, where $F\in C^2$, $X\in\RR^n$ for some $n\in\mathbb{N}$. Assume that $X_0\in\RR^n$  is an equilibrium  and $\gamma$ is a non-stationary solution of the system.
            The trajectory/orbit $\gamma$ is said to be \emph{\textbf{homoclinic}} to $X_0$ if $$\gamma\left(t\right)\to X_0 \quad\text{ as }\quad t\to\pm\infty.$$
        \end{definition}
        The existence of a saddle point equilibrium is not a sufficient condition for the existence of a homoclinic trajectory, as this also requires the existence of an intersection between the stable and unstable manifolds of the equilibrium point.\\
       
\begin{figure}
            \centering
            \begin{subfigure}{.5\linewidth}
                \centering
                \includegraphics[width=\linewidth]{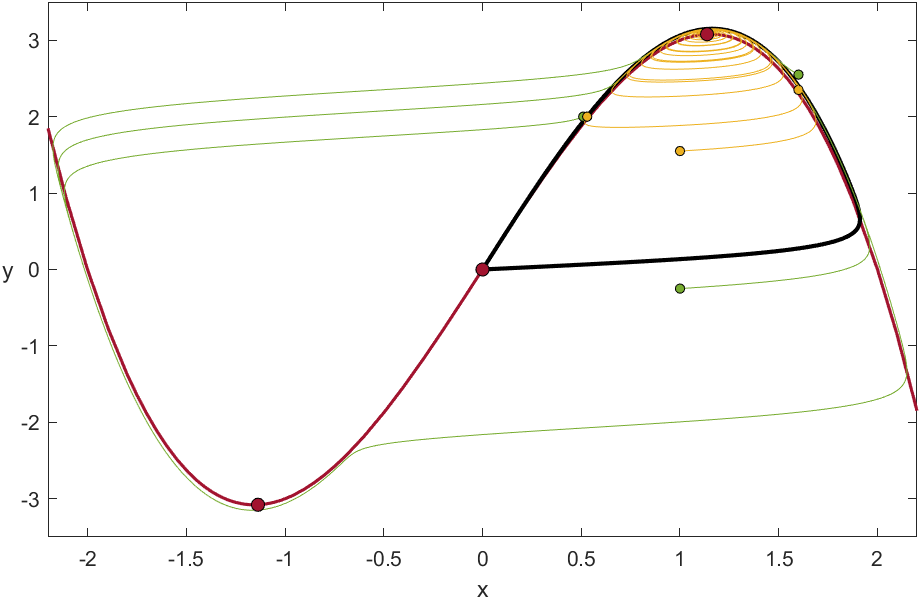}
                \caption{ \small }
                \label{fig:HomOrb_xy}
            \end{subfigure}\hfill
            \begin{subfigure}{.41\textwidth}
                \centering
                \includegraphics[width=\linewidth]{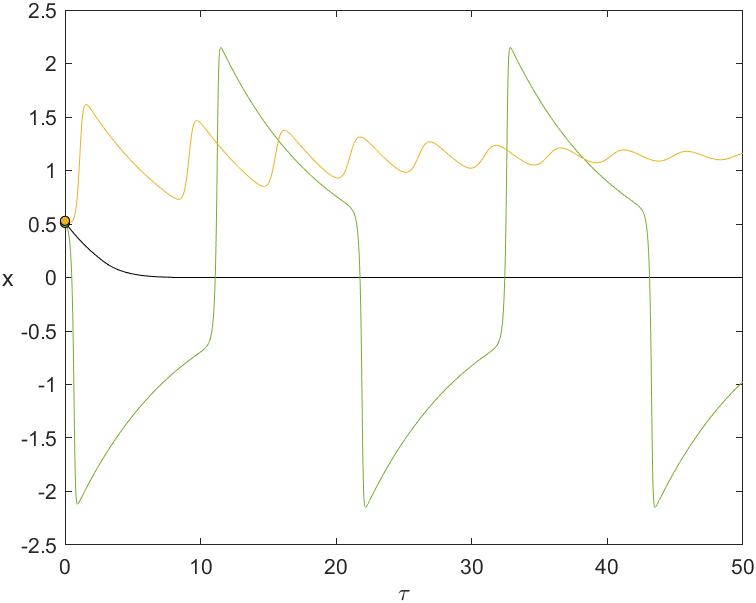}
                \caption{ \small }
                \label{fig:HomOrb_xt}
            \end{subfigure}
            \caption{ \small  (a) - $\left(x,y\right)$ plane of the system \eqref{fhn_c=0} for $\varepsilon=0.5$ and $b=b_{0.5}^{Hom}\approx 0.36932$. $\VC$ and equilibrium points in red. Homoclinic orbit (black), trajectories whose initial condition is within the region bounded by the homoclinic curve (yellow) and trajectories whose initial condition is outside the region bounded by the homoclinic curve (green). (b) - $\left(\tau,x\right)$ plane of the system \eqref{fhn_c=0} for $\varepsilon=0.5$ and $b=b_{0.5}^{Hom}\approx 0.36932$. Trajectories initiated at $\left(x_0,y_0\right)=\left(0.51,2\right)$ (green), at $\left(x_0,y_0\right)=\left(0.53,2\right)$ (yellow) and at the homoclinic orbit $\left(x_0,y_0\right)\approx\left(0.5158,1.9578\right)$.}
            \label{fig:HomOrbs}
        \end{figure}

        We already knew that \eqref{fhn_c=0} has the saddle equilibrium $\mathrm{E}_0=(0,0)$ and we have observed numerically that, for $b=b_\varepsilon^{Hom}$, there is a homoclinic trajectory to the origin $\mathrm{E}_0$. In Figures \ref{fig:bifb_eps05_mais_exemplo} and \ref{fig:HomOrbs}, which consider the system \eqref{fhn_c=0} for $\varepsilon=0.5$, the homoclinic orbit appears for $b=b_\varepsilon^{Hom}\approx 0.36932$, very close to $b_{\varepsilon}^H\approx 0.3666$.
        The bifurcation where a periodic solution is destroyed at a homoclinic orbit is called a \emph{\textbf{homoclinic bifurcation}} \cite{guk,strogatz}.

    %
\section{ Singularities and {Canards}}\label{sec:canards}
    As we have seen in case (i) of Section \ref{sec:bif}, when the system's equilibrium is also a fold point, a subcritical Hopf bifurcation occurs, giving rise to a stable periodic solution (cf. Figure \ref{fig:bifc_eps05_mais_exemplo}).
    In Sections \ref{sec:sing} and \ref{sec:se}, we found another periodic solution which  is not related to the Hopf bifurcation.
    Thus, for different parameters, we have two stable periodic solutions with distinct dynamical properties.
    In this section, the \emph{\bf canard phenomenon} bridges the local periodic solution from the Hopf bifurcation to the global limit cycle.
    %
    
    As illustrated in Figure \ref{fig:b0_c1152_e05}, the existence in \eqref{fhn_b=0} of a periodic solution around the equilibrium point implies that the orbit traverses a part of its trajectory close to the repelling region of the critical manifold $\VC$. This behaviour seems to be contradictory since we found a solution that is close to the unstable branch of $\VC$.
    \begin{definition}
        A slow trajectory that is $\mathcal{O}\left(\varepsilon\right)$ away from the repelling region of the slow manifold for a time of order $\mathcal{O}\left(1\right)$, is called a \emph{\textbf{canard}}.
    \end{definition}

    \emph{Canards} are related to equilibria of the slow equation that occur at fold points of the slow manifold.
    \begin{definition}\label{def:fold_point}
        Given a singular $(1,1)$--fast-slow system parametrised by $\lambda\in\RR$. Let $\p=\left(x,y\right)$ be a fold point. The point $\p$ is called a \textbf{singular fold} if:
        \begin{equation}
            \label{dobrasingular}
            \begin{aligned}
                &f\left(\p,\lambda,0\right)=0, && f_x\left(\p,\lambda,0\right)=0,\\
                &f_{xx}\left(\p,\lambda,0\right)\neq0, && f_y\left(\p,\lambda,0\right)\neq0 \quad \text{and} \quad g\left(\p,\lambda,0\right)=0.
            \end{aligned}
        \end{equation}     
    \end{definition}

    The definition of \emph{singular fold point} differs from that of a fold point since it   is required to coincide with an equilibrium point of system \eqref{deflento} (i.e., $g\left(\p,\lambda,0\right)=0$).
    \begin{definition}\label{def:regular_fold}
        A singular fold $\p$ is said to be \textbf{regular} if:
        \begin{equation}
            g_x\left(\p,\lambda,0\right)\neq0 \quad \text{ and } \quad g_{\lambda}\left(\p,\lambda,0\right)\neq0.
        \end{equation}
    \end{definition}

    The behaviour around a regular singular fold is described in the next theorem, a concatenation of two important results, whose proof can be found in \cite{KS01a,KS01b,kuehn}. 
        
     \begin{theorem}[Krupa \& Szmolyan, 2001a, 2001b \cite{KS01a,KS01b}.
    See also Theorems 8.1.3 and 8.2.1 of 
    \citeyear{kuehn}]\label{th:canard}
        Consider a (1,1)--fast-slow system where $\left(x,y\right)=\left(0,0\right)$ is a regular singular fold point for $\lambda=0$, in the following normal form:
        \begin{equation}\label{eq:nformHopfsing}
            \begin{aligned}
                &x'=-y\,l_1\left(x,y,\lambda,\varepsilon\right)+x^2\,l_2\left(x,y,\lambda,\varepsilon\right) + \varepsilon\,l_3\left(x,y,\lambda,\varepsilon\right)\\
                &y'=\varepsilon \left(\pm x\,l_4\left(x,y,\lambda,\varepsilon\right)-\lambda \,l_5\left(x,y,\lambda,\varepsilon\right)+y\,l_6\left(x,y,\lambda,\varepsilon\right)\right).
            \end{aligned}
        \end{equation}
        %
    
        Assume that, for $\varepsilon=0$, there is a slow trajectory connecting the repelling and attracting regions of the critical manifold $\VC$.
        Then, there exist $\varepsilon_0>0$ and $\lambda_0>0$ such that for $0<\varepsilon<\varepsilon_0$ and $|\lambda|<\lambda_0$, the system has an equilibrium point $\p\in\RR^2$ near the origin where $\p\to\left(0,0\right)$ as $\left(\lambda,\varepsilon\right)\to\left(0,0\right)$.
    
        Moreover, if $\p$ is stable for $\lambda<0$, there exists a continuous function $\lambda_c:\left[0,\varepsilon_0\right]\to\RR$ that associates each value of $\varepsilon\in \left[0,\varepsilon_0\right]$ to a value $\lambda$ that gives rise to a \emph{canard} in the system, asymptotically defined by: 
        \begin{equation*}
            \lambda_c\left(\sqrt{\varepsilon}\right)=-\left(B + A\right)\varepsilon+\mathcal{O}\left(\varepsilon^{3/2}\right),    
        \end{equation*}
        and there exists a continuous function $\lambda_\mathrm{H}:\left[0,\varepsilon_0\right]\to\RR$ that associates each value of $\varepsilon\in \left[0,\varepsilon_0\right]$ to a value $\lambda$ that gives rise to a Hopf bifurcation in the system, asymptotically defined by: 
        \begin{equation*}
            \lambda_\mathrm{H}\left(\sqrt{\varepsilon}\right)=-B\varepsilon+\mathcal{O}\left(\varepsilon^{3/2}\right),    
        \end{equation*}
        where
        \begin{equation*}
            A=\dfrac{-\dfrac{\partial l_1}{\partial x}+3\dfrac{\partial l_2}{\partial x}-2\dfrac{\partial l_4}{\partial x}+2l_6}{8}\qquad B=\dfrac{\dfrac{\partial l_3}{\partial x}+l_6}{2},    
        \end{equation*}
        where the functions $l_i,\,i=1,\ldots,6$ and their partial derivatives are evaluated at the point $\left(x,y,\lambda,\varepsilon\right)=\left(0,0,0,0\right)$.
         
        For $A\neq 0$, the Hopf bifurcation is non-degenerate. The Hopf bifurcation is supercritical if $A<0$, and subcritical if $A>0$. 
    \end{theorem}
     
    Whenever a system with a regular singular fold exists, there is a change of coordinates that transforms the system into \eqref{eq:nformHopfsing}. 
    As noted in Section \ref{th:canard}, the existence of \emph{canards} is coupled with the existence of a Hopf bifurcation, which in such systems is referred to as a \textbf{singular Hopf bifurcation}.
    That is, the values of $\lambda_c$ that give rise to \emph{canards} are at most  $\mathcal{O}\left(\varepsilon\right)$ away from the values $\lambda_\mathrm{H}$ where a Hopf bifurcation occurs,
    \begin{equation*}
        \lambda_\mathrm{H}-\lambda_c=A\,\varepsilon+\mathcal{O}\left(\varepsilon^{3/2}\right).    
    \end{equation*}
    An important consideration is the empirical difficulty in finding \emph{canards}, since the smaller $\varepsilon$ is, the narrower the interval of $\lambda$ values for which \emph{canards} appear.
    This small interval separates the $\lambda$ values at which a limit cycle occurs  from those where it either undergoes a Hopf bifurcation or does not have any periodic solution.
    { Figure \ref{fig:bifc_eps05_diagram} and its magnified Figures \ref{fig:bifc_eps05_diagram_amp} and \ref{fig:bifc_eps05_diagram_2xamp} illustrate this phenomenon.}
    For this reason, the dynamics produced in the system around this small interval is called the \emph{\textbf{canard explosion}}.
    \begin{example}\label{ex:super_Hopf}
      Consider the case (i) of Section \ref{sec:bif}. Since at the Hopf bifurcation point $c^\mathrm{H}$ the equilibrium $\mathrm{E}=\left(c^\mathrm{H},\varphi\left(c^\mathrm{H}\right)\right)=\left(\dfrac{2}{\sqrt{3}},\dfrac{16}{3\sqrt{3}}\right)$ is a regular fold singularity, to determine the stability of the periodic solution, we need to transform the system into its normal form (see \citet{guk, kuehn} for more details), which involves a translation of the equilibrium and of the bifurcation parameter to the origin. More specifically, we consider the change of coordinates:
        \begin{equation*}
            \left(\bar{x},\bar{y},\bar{c}\right) = \left(x-c^\mathrm{H},y-\varphi\left(c^\mathrm{H}\right), c^\mathrm{H}-c\right),
        \end{equation*}
        where $\varphi\left(c^\mathrm{H}\right)=4c^\mathrm{H}-\left(c^\mathrm{H}\right)^3$. Omitting the bars $\bar{x}$ and $\bar{y}$ this gives rise to the following system equivalent to \eqref{fhn_b=0}:
        \begin{equation}
            \begin{aligned}
                &x'=-y+x^2\left(-2\sqrt{3}-x\right)\\
                &y'=\varepsilon \left(x+\bar{c}\right).
            \end{aligned}
        \end{equation}

        Since the parameters $c$ and $\bar{c}$ have opposite orientations, then the existence of a supercritical Hopf bifurcation with periodic solutions for $\bar{c}>0$ implies the existence of a subcritical Hopf bifurcation with periodic solutions for $c<c^H$ in the original system \eqref{fhn_b=0}.
            
        It is easy to verify that $\p=\left(-\bar{c},\varphi\left(-\bar{c}\right)\right)\,\overset{\bar{c}\to 0}{\longrightarrow}\left(0,0\right)\,$ is an equilibrium of the system and that $\p$ is stable for $\bar{c}<0$. Additionally, we have
        \begin{equation*}
            f_x\left(\p,0\right)=0,\quad f_{xx}\left(\p,0\right)=-4\sqrt{3}\neq0\quad \text{and}\quad  f_y\left(\p,0\right)=-1\neq0,    
        \end{equation*}
        and also,
        \begin{equation*}
            g_x\left( \p,0\right)=1\neq0 \quad \text{and}\quad g_{\Bar{c}} \left(\p,0\right)=1\neq0,    
        \end{equation*}
        so $\p$ is a regular singular fold point for $\Bar{c}=0$. 
        Note that the system is defined in the normal form required by Theorem \ref{th:canard} and we have:
        \begin{equation*}
            l_1=1,\quad l_2=-2\sqrt{3}-x,\quad l_3=0,\quad l_4=1,\quad l_5=-1,\quad l_6=0\quad \text{and}\quad A=-3/8<0.    
        \end{equation*}
        Thus, the Hopf bifurcation is supercritical for $\Bar{c}$ (subcritical for $c$) and occurs when 
        \begin{equation*}
            \bar{\mathrm{c}}_\mathrm{H}\left(\sqrt{\varepsilon}\right)=0\cdot\,\varepsilon + \mathcal{O}\left(\varepsilon^{3/2}\right)    
        \end{equation*} \noindent and the \emph{canards} occur for 
        \begin{equation*}
            \bar{\mathrm{c}}_c\left(\sqrt{\varepsilon}\right)=-\dfrac{3}{8}\,\varepsilon + \mathcal{O}\left(\varepsilon^{3/2}\right).\qquad\diamond    
        \end{equation*}
    \end{example}

        \begin{example}\label{ex:Hopf_c0}
            Now consider the case (ii) of Section \ref{sec:bif} (with $c=0$) and  the equilibrium $E_+=(x_+,x_+/b)$, with $x_+=\sqrt{4-b^{-1}}$.
            Then $E_+$ is a regular singular fold point for $b=3/8$ with $x_+=4/3$ and there is a Hopf bifurcation point at $b=b_\varepsilon^H$ with $b_\varepsilon^H\to 3/8$ as $\varepsilon\to 0$.
            In order to apply Theorem~\ref{th:canard} we change variables and parameters by
            $$
            (\bar{x},\bar{y},\lambda)=\left(x-x_+,y-8x_+/3,b-3/8 \right) 
            $$
            to obtain the equivalent system
            \begin{equation}\label{eq:nfEx62}
                \begin{array}{lcl}
                    \bar{x}'&=&-\bar{y}-\bar{x}^3-3\bar{x}^2x_+\\
                    \bar{y}'&=&\varepsilon\left(\bar{x}-\lambda\bar{y}-\frac{3}{8}\bar{y} \right) .
                \end{array}
            \end{equation}
            In the notation of  Theorem~\ref{th:canard} we have
            $$
            l_1=-1\qquad
            l_2=-3x_+-\bar{x}\qquad
            l_3=0\qquad
            l_4=1\qquad
            l_5=1\qquad
            l_6=-\frac{3}{8}
            $$
            hence $A=-\dfrac{1}{8}\left(3+\dfrac{3}{4}\right)<0$ therefore the bifurcation is supercritical for $\lambda$, and also for $b$ in the original system \eqref{fhn_c=0}.
            Since for $b<b_\varepsilon^H$ the eigenvalues of the matrix $J_{E_+}$ are positive, then 
            the bifurcating periodic solution is unstable.
        \end{example}

\begin{figure}[h]
        \centering
        \begin{subfigure}{.45\linewidth}
            \centering
            \includegraphics[width=\linewidth]{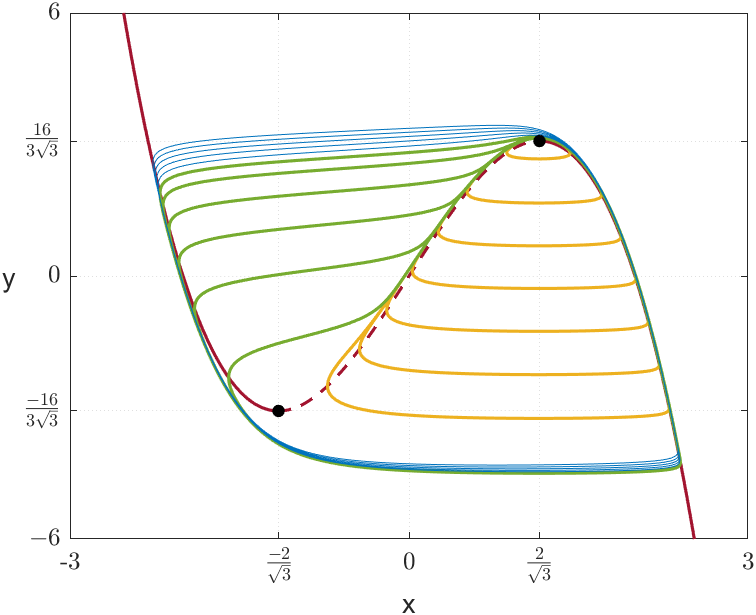}
            \caption{ \small }
            \label{fig:bifc_eps05_canards_2d}
        \end{subfigure}%
        \hfill
        \begin{subfigure}{.45\textwidth}
            \centering
            \includegraphics[width=\linewidth]{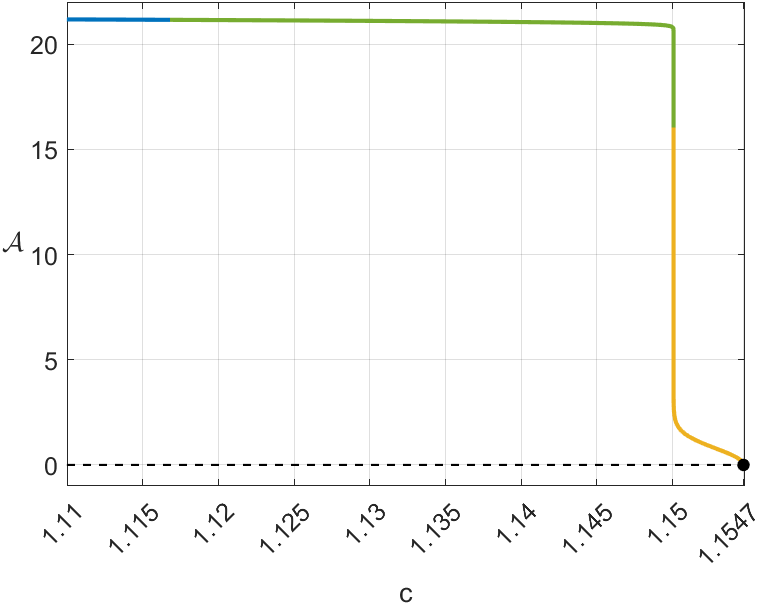}
            \caption{ \small }
            \label{fig:bifc_eps05_diagram_canard}
        \end{subfigure}
        \vfill
        \begin{subfigure}{.45\linewidth}
            \centering
            \includegraphics[width=\linewidth]{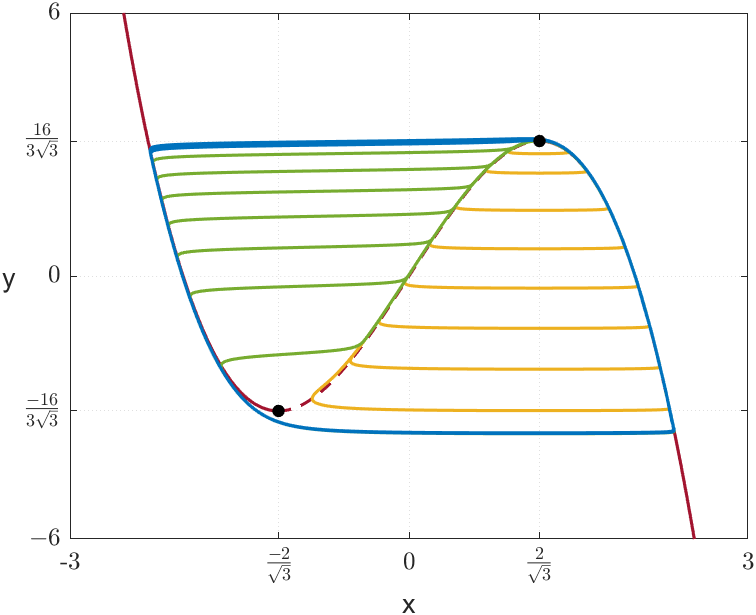}
            \caption{ \small }
            \label{fig:bifc_eps01_canards_2d}
        \end{subfigure}%
        \hfill
        \begin{subfigure}{.45\textwidth}
            \centering
            \includegraphics[width=\linewidth]{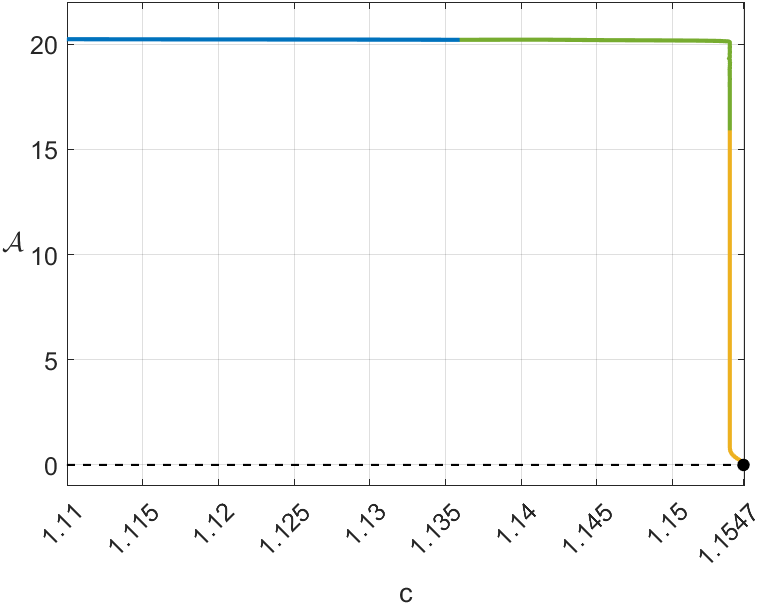}
            \caption{ \small }
            \label{fig:bifc_eps01_diagram_canard}
        \end{subfigure}
        \caption{ \small  (a) - {Periodic solutions} of the system \eqref{fhn_b=0} {for several values of $c$}, with $\varepsilon=0.5$. Critical manifold in red, attracting regions (solid line) and repelling region (dashed line), fold points (black), \emph{canards} without head (yellow), \emph{canards} with head (green), and limit cycles (blue). (b) -  Evolution of the   length  $\mathcal{A}$ of the \emph{canards} of (a) as function of $c$.   Unstable equilibrium point (dashed line), Hopf bifurcation point (black dot),  length of \emph{canards} without head and Hopf orbits (yellow),  length of \emph{canards} with head orbits (green), and   length of limit cycles (blue). (c) - { Periodic solutions} of the system \eqref{fhn_b=0} {for several values of $c$}, with §$\varepsilon=0.1$ (conventions as in (a)) (d) -  {Evolution of the  length of the \emph{canards} of (c) as a function of the parameter $c$, conventions as in (b).}}  
        \label{fig:canards}
    \end{figure}

    Figure \ref{fig:canards} illustrates the emergence of \emph{canards} in system \eqref{fhn_b=0} with $\varepsilon=0.5$ and $\varepsilon=0.1$.
    As we may see, the stable periodic solutions emerge from the subcritical Hopf bifurcation for $c=2/\sqrt{3}$ (analogous for $c=-2/\sqrt{3}$).
    As mentioned above, some of these periodic solutions remain close to the repelling region of the critical manifold, making them \emph{canards}.
    In Figures \ref{fig:bifc_eps05_canards_2d} and \ref{fig:bifc_eps01_canards_2d}, two types of orbits can be distinguished: the first type, coloured in yellow, consists of slow trajectories near $\mathcal{C}_{0\text{M}}$ and $\mathcal{C}_{0\text{R}}$ and one fast trajectory;
    the second type, coloured in green, which \emph{surrounds} the other fold point, consists of three slow trajectories near $\mathcal{C}_{0\text{L}}$, $\mathcal{C}_{0\text{M}}$ and $\mathcal{C}_{0\text{R}}$ and two fast trajectories.
    We refer to these two types of \emph{canards} as \textbf{headless \emph{canard}} and \textbf{headed \emph{canard}}, respectively\footnote{This terminology is motivated by the literal meaning of the word ``canard" in French (\emph{canard} = duck). In this system, it is a bit more complicated to understand why, but if the reader tilts their head to the left they will see that the \emph{canard} orbits resemble the shape of a duck, with or without head.}.
         
    Note that the  length of the \emph{canards} grows almost instantaneously.
    As stated in \citet{kuehn}, this was expected, but it is still surprising.
    The \textit{nearly vertical} lines observed in Figures \ref{fig:bifc_eps05_diagram_canard} and \ref{fig:bifc_eps01_diagram_canard} show \emph{canards} with  lengths, $\mathcal{A}\left(c\right)\in\left[2,20\right]$ for $c\approx 1.150077$ in the case $\varepsilon=0.5$ and for $c\approx 1.153794$ in the case $\varepsilon=0.1$.
    %
\section{Discussion and Further Work}\label{discussion}
    This article   studies the dynamics of a fast-slow system derived from the FH-N model according to the geometric singular perturbation theory and the bifurcation theory methods. In particular, we  provide an analytic proof that the fast-slow FH-N system presents relaxation oscillation dynamics as well as periodic solutions induced by Hopf bifurcation. We   emphasise the emergence of a homoclinic orbit and \emph{canards} connecting these two distinct dynamics. We  illustrate each result with numerical computations.
    In each section, we have complemented the discussion of observed dynamics with stability and bifurcation analysis. This approach empowers the readers to navigate the parameter space, enabling them to find specific dynamical regimes of interest.

    The FH-N equation is an excellent candidate to serve as a starting point for developing new methods to analyse transient and non-reciprocal coupled interactions. Cardiac and neuronal examples, in particular, highlight the FH-N model's adaptability and significance in these areas \cite{lacasa24}.
    
    When coupling two FH-N systems through the slow equations, we have performed numerical simulations that suggest the existence of mixed-mode oscillations where the dynamics is characterised by the existence of orbits that alternate periodically between large and small amplitude oscillations \cite{GLR24}. 
    Recently, \citet{FastCoupling_Kristian} showed that the coupling through the fast equations of two distinct FH-N models produces mixed-mode oscillations induced by a cusp singularity present in this system. We believe that a similar result may explain our numerical findings.  

    The prevalence of chaotic attractors in periodically perturbed fast-slow systems is yet to be explored. The article \citet{GuckYoung} bridged the areas of geometric singular perturbation theory and of chaotic attractors theory and provided a general technique for proving the existence of chaotic attractors for periodically perturbed two-dimensional vector fields with two time scales, that are equivalent to:
    \begin{equation}\label{formula_discussion}
        \begin{aligned}
            &\varepsilon\,\dfrac{\mathrm{d}x}{\mathrm{d}t}= f(x,y,\theta)\\
            &\phantom{\varepsilon\,}\dfrac{\mathrm{d}y}{\mathrm{d}t} = g(x,y,\theta)\\
            &\phantom{\varepsilon\,}\dfrac{\mathrm{d}\theta}{\mathrm{d}t}= \omega,\\
        \end{aligned}
    \end{equation}
    where $(x,y, \theta)\in \RR\times \RR\times \EU^1$ and $\omega >0$ is the slow driving frequency. 
    
    Results in \citet{GuckYoung} connect two areas of dynamical systems: the theory of chaotic attractors for discrete two-dimensional H\'enon-like maps and geometric singular perturbation theory. Two-dimensional H\'enon-like maps are diffeomorphisms whose singular limit is a family of non-invertible one-dimensional maps. In \citet{Wang_Young}, the authors obtained sufficient conditions for the existence of chaotic attractors in these families. Three-dimensional singularly perturbed vector fields have return maps that are also two-dimensional diffeomorphisms with folds. To represent fully the behaviour of system \eqref{formula_discussion} in the singular limit, we must allow \emph{jumps} of trajectories from one sheet of the critical manifold to another that follow the direction of trajectories when $\varepsilon> 0$. \emph{Jumps} parallel to the $x$-axis occur at folds, where the tangent plane to the critical manifold includes this direction. We have chosen to leave the study of periodically perturbed FH-N models for future work.

\section*{Acknowledgments}
The  first and second authors were partially supported by CMUP, member of LASI, which is financed by national funds through FCT – Fundação para a Ciência e a Tecnologia, I.P., under the projects with reference UIDB/00144/2020 and UIDP/00144/2020. 
The third author has been supported by the Project CEMAPRE/REM – UIDB/05069/2020 financed by FCT/MCTES through national funds. 
    
\end{document}